\documentclass[10pt]{article}
\usepackage[top=1in, bottom=1in, left=1in, right=1in]{geometry}
\usepackage{srcltx,graphicx}
\usepackage{amsmath, amssymb, amsthm}
\usepackage{color, xcolor, colortbl}
\usepackage{subfig, overpic}
\usepackage{hyperref}
\usepackage{epstopdf}
\usepackage{algorithm}
\usepackage{algpseudocode}
\usepackage[capitalize,nameinlink]{cleveref} %% use ``cref'' to replace ``ref

\usepackage{stackengine}
\usepackage{scalerel}
\newcommand\trianglesignplus[1][2.3ex]{%
    \renewcommand\stacktype{L}%
    \mathbin{%
    \scaleto{\stackon[1.2pt]{$\triangle$}{\tiny $\boldsymbol{+}$}}{#1}%
}}

\graphicspath{{images/}}
\numberwithin{equation}{section}

\newcommand\bbR{\mathbb{R}}
\newcommand\bbZ{\mathbb{Z}}

\newcommand\cM{\mathcal{M}}
\newcommand\cN{\mathcal{N}}
\newcommand\cT{\mathcal{T}}

\newcommand\dd{\,\mathrm{d}} %differential operator

\newcommand\sps[1]{^{(#1)}}
\newcommand\wl[2]{\sps{#1}_{#2}}

\newcommand\Ntrainsample{{N^{\mathrm{train}}_{\mathrm{samples}}}}
\newcommand\Ntestsample{{N^{\mathrm{test}}_{\mathrm{samples}}}}
\newcommand\Nparams{{N_{\mathrm{params}}}}
\newcommand\trainerror{{\epsilon_{\mathrm{train}}}}
\newcommand\testerror{{\epsilon_{\mathrm{test}}}}
\newcommand\ipt{{\mathrm{in}}}
\newcommand\out{{\mathrm{out}}}

\newcommand\sffont[1]{{\sf{#1}}}

\newcommand\id{{\sffont{id}}}
\newcommand\Conv{{\sffont{Conv }}}
\newcommand\LC{{\sffont{LC }}}
\newcommand\Dense{{\sffont{Dense }}}

\newcommand\eps{\varepsilon}
\newcommand\argmin{\mathrm{argmin}\,}
\newcommand\NN{\mathrm{NN}}

\title{ BCR-Net: a neural network based on the nonstandard wavelet form }

\date{}
\author{
  Yuwei Fan\thanks{Department of Mathematics, Stanford University,
    Stanford, CA 94305. Email: {\tt ywfan@stanford.edu}},~~
  Cindy Orozco Bohorquez\thanks{Institute for Computational and Mathematical Engineering, 
    Stanford University, Stanford, CA 94305. Email: {\tt orozcocc@stanford.edu}}, ~~
  Lexing Ying\thanks{Department of Mathematics and ICME, Stanford University,
    Stanford, CA 94305.
    Facebook AI Research, Menlo Park, CA 94025.
    Email: {\tt lexing@stanford.edu}}
}

\begin{document}
\maketitle

\begin{abstract}
  
  This paper proposes a novel neural network architecture inspired by the nonstandard form proposed
  by Beylkin, Coifman, and Rokhlin in [Communications on Pure and Applied Mathematics, 44(2),
    141-183]. The nonstandard form is a highly effective wavelet-based compression scheme for linear
  integral operators. In this work, we first represent the matrix-vector product algorithm of the
  nonstandard form as a linear neural network where every scale of the multiresolution computation
  is carried out by a locally connected linear sub-network. In order to address nonlinear problems,
  we propose an extension, called BCR-Net, by replacing each linear sub-network with a deeper and
  more powerful nonlinear one. Numerical results demonstrate the efficiency of the new architecture
  by approximating nonlinear maps that arise in homogenization theory and stochastic computation.
  
  %By taking advantage of the
%multiscale structure of BCR algorithm, we construct a hierarchical network consisting of a group of
%sub-networks, each of which corresponds to a scale of wavelet transform, and is a locally connected
%network.
%(here ``locally connected'' doesn't equal to the LC defined in the paper. The point is ``local'', rather than ``globally''
%where each level reproduces a wavelet filter by combining single layer convolutional and locally
%connected networks. 

%\vspace*{4mm}
\end{abstract}


%==================================================================================================
\section{Introduction}\label{sec:intro}

%In this paper, we introduce a novel neural network architecture inspired by the wavelet-based fast
%algorithms.

This paper proposes a novel neural network architecture based on wavelet-based compression
schemes. There has been a long history of research on representing differential and integral
operators using wavelet bases. Such a representation is particularly attractive for integral
operators such as pseudo-differential operators and Calderon-Zygmund operators because the
wavelet-based representation is often sparse due to the vanishing moment property of the wavelet
basis \cite{cohen2003numerical,daubechies1988orthonormal}. For a typical problem discretized with
$N$ unknowns, a direct application of the wavelet transform, followed by a thresholding step,
results in $O(N\log N)$ significant matrix entries, with the prefactor constant depending
logarithmically on the accuracy level.

%% particularly the algorithm proposed in 1991 by Beylkin,
%% Coifman and Rokhlin, and nowadays known as the BCR algorithm \cite{bcr}. Before the development of
%% BCR, wavelet bases were already extensively used to compress dense linear operators such as
%% pseudo-differential operators and Calderon-Zygmund operators. Nevertheless, 
%% the vast majority of the existing algorithms provided either a
%% fast algorithm for a specific operator \cite{alpert1991fast, carrier1988fast, greengard1987fast},

In a ground-breaking paper \cite{bcr}, Beylkin, Coifman, and Rokhlin proposed a nonstandard form of
the wavelet representation that surprisingly reduced the number of significant entries to
$O(N)$. The key idea of this work is rather simple: instead of treating the matrix obtained from
integral operator discretization as a linear map, view it as a two-dimensional image and hence
compress it using two-dimensional wavelets. Not only drastically reducing the number of significant
entries, but the resulting algorithm for matrix-vector multiplication is also significantly easier to
implement and applicable to rather general pseudo-differential operators. One natural but rather
unexplored question is how to extend the nonstandard wavelet form to nonlinear integral operators
(see \cite{alpert2002adaptive,yavneh2006multilevel} for some related work).

Recently neural networks (NNs) have experienced great successes in artificial intelligence
\cite{Hinton2012, Krizhevsky2012, goodfellow2016deep, MaSheridan2015, Leung2014, SutskeverNIPS2014,
  leCunn2015, SCHMIDHUBER2015} and even in solving PDEs \cite{khoo2017solving, berg2017unified,
  han2018solving, fan2018mnn, fan2018mnnh2, Araya-Polo2018, Raissi2018}. One of the key reasons for
these successes is that deep neural networks offer a powerful way for approximating high-dimensional
functions and maps. Given a nonlinear integral operator, a fully connected NN is theoretically
capable of approximating such a map \cite{CohenSharir2018,Hornik91,Khrulkov2018,Mhaskar2018}.
Nonetheless, using a fully connected NN often leads to a prohibitively large number of parameters,
hence long training stages and overwhelming memory footprints. To overcome these challenges, one can
incorporate knowledge of the underlying structure of the problem in designing suitable network
architectures. One promising and general strategy is to build NNs based on a multiscale
decomposition~\cite{fan2018mnnh2,fan2018mnn,Yingzhou2018,WangChung2018}. The main idea, frequently
used in image processing as well
\cite{badrinarayanan2015segnet,Bruna2012,Chen2018DeepLab,LITJENS201760,Ronneberger2015,Ulyanov2018},
is to learn increasingly coarse-grained features of a complex problem across different layers of the
network structure, so that the number of parameters in each layer can be effectively controlled.

In this work, we introduce a novel NN architecture based on the nonstandard wavelet form in
\cite{bcr} in order to approximate certain global nonlinear operators. The paper is organized as
follows. \cref{sec:non} describes the nonstandard wavelet form and the fast matrix-vector
multiplication based on it.  \cref{sec:mvnn} introduces some basic NN tools and represent the
matrix-vector multiplication algorithm in the form of a linear NN. \cref{sec:bcrnet} generalizes the
linear NN to the nonlinear case by incorporating nonlinear activation functions and increasing the
number of layers and channels of the NN. \cref{sec:application} describes implementation details and
demonstrates the numerical efficiency of the NN with two applications from homogenization theory and
stochastic computation. In both cases, the nonlinear mapping can be well approximated by the
proposed NN, at a relative accuracy $10^{-4}\sim 10^{-3}$, with only $10^4\sim 10^{5}$ parameters
for 2D problems and $10^5$ parameters for 3D problems.

%We demonstrate the capabilities of MNN-BCR with two classical yet challenging examples in
%computational physics:\dots.

%%   Beylkin-Coifman-Rokhlin (BCR) algorithm
%% \cite{bcr} is a wavelet-based algorithm for the matrix-vector multiplication, in which the matrix
%% (whenever regularity permits) is converted to a sparse form in the wavelet bases in the
%% \emph{no-standard form}. For pseudo-differential operators and Calderon-Zygmund operators, the
%% matrix can be compressed to a sparse matrix with only $O(n)$ nonzero entries.

%thus the arithmetic complexity of the matrix-vector multiplication is reduced to $O(n)$ with the
%help of fast wavelet transform \cite{mallat1989theory}.

%% to effectively approximate nonlinear
%% maps in the form
%% \begin{equation}\label{eq:map}
%%     u = \cM(v), \quad u, v\in\Omega\subset\bbR^{n},
%% \end{equation}
%% where $u$ is the quantity of interest and $v$ is the parameter that serves to identity a particular
%% configuration of the system. This map is ubiquitous from the discretization of integral equations,
%% parameterized partial differential equations and pseudo-differential operators. 
%% We focus on the fast and efficient evaluation of the map \eqref{eq:map}.

%==================================================================================================
\section{Nonstandard wavelet form}\label{sec:non}

%In this section, we aim to depict the matrix-vector multiplication of BCR algorithm within the
%framework of neural networks.

%We briefly review the BCR algorithm for the 1D case in \cref{sec:bcr} and introduce a linear neural
%network architecture for the matrix-vector multiplication of BCR algorithm in \cref{sec:matvec}. An
%extension to the multi-dimensional setting is presented in \cref{sec:nD}.

In this section, we briefly summarize the nonstandard wavelet form proposed in \cite{bcr}, using the
compactly supported orthonormal Daubechies wavelets (see \cite{daubechies1988orthonormal} for
details) as the basis functions.

%============
\subsection{Wavelet transform}
In a one-dimensional multiresolution analysis, the starting point is a {\em scaling function}, or
sometimes called father wavelet, $\varphi(x)$ that generates a family
\begin{equation}\label{eq:translate_phi}
  \varphi_{k}^{(\ell)}(x)=2^{\ell/2}\varphi(2^{\ell}x-k),\quad
  \ell=0,1,2,\dots,\quad k\in\bbZ,
\end{equation}
such that for each fixed $\ell$ the functions $\{\varphi^{(\ell)}_k\}_{k\in\mathbb{Z}}$ form a Riesz
basis for a space $V_{\ell}$. The spaces $\{V_{\ell}\}_{\ell \ge 0}$ form a nested sequence of
spaces $V_{\ell}\subset V_{\ell+1}$. Because of this nested condition, $\varphi(x)$ satisfies the
\textit{dilation relation}, also known as the \textit{refinement equation}
\begin{equation}\label{eq:recursive_phi}
  %\frac{1}{\sqrt{2}} \varphi\left(\frac{x}{2}\right) = \sum_{i\in\bbZ} h_i\varphi(x-i).
  \varphi(x) = \sqrt{2} \sum_{i\in\bbZ} h_i\varphi(2x-i).
\end{equation}
In the case of Daubechies wavelets, $\varphi(x)$ is supported in $[0,2p-1]$ for a positive integer
$p$ and hence the coefficients $\{h_i\}$ are nonzero only for $i=0,\ldots, 2p-1$. From the
orthonormal condition, the function $\varphi(x)$ and its integer translates satisfy an orthonormal
condition
\begin{equation}\label{eq:normalized}
  \int_{\bbR}\varphi(x-a)\varphi(x-b)\dd x = \delta_{a,b}.
\end{equation}
In terms of the coefficients $h_i$ in \eqref{eq:recursive_phi}, this orthonormal condition can be
written as
\begin{equation}\label{eq:orthogonality}
  \sum_{i\in\bbZ} h_i^2=1,\quad 
  \sum_{i\in\bbZ} h_ih_{i+2m}=0, \quad m\in\bbZ \backslash \{0\},
\end{equation}
following \cite{daubechies1988orthonormal}.

Given the scaling function $\varphi(x)$, the \emph{mother wavelet function} $\psi(x)$ is defined as
\begin{equation}\label{eq:recursive_psi}
  %\frac{1}{\sqrt{2}} \psi\left(\frac{x}{2}\right) = \sum_{i\in\bbZ} g_i\varphi(x-i),
  \psi(x) = \sqrt{2} \sum_{i\in\bbZ} g_i\varphi(2x-i),
\end{equation}
with $g_i=(-1)^{1-i} h_{1-i}$ for $i\in\bbZ$. From the support of $\varphi(x)$ and the non-zero
pattern of $h_i$, it is clear the support of $\psi(x)$ is $[-p+1,p]$ and $g_i$ is nonzero only for
$i=-2p+2,\ldots,1$.  The {\em Daubechies wavelets} are defined as the scaled and shifted copies of
$\psi(x)$

%% To clarify the notations introduced above, we take Haar wavelet as an example.
%% \begin{example}[Haar wavelet]
%% For Haar wavelet, the scaling function and wavelet function are
%% \begin{equation}
%%     \varphi(x)= \begin{cases}
%%         1,  & x\in [0,1],\\
%%         0,  & \text{elsewhere},
%%     \end{cases}
%%     \quad
%%     \psi(x)= \begin{cases}
%%         1,  & x\in [0,1/2],\\
%%         -1,  & x\in (1/2,1],\\
%%         0,  & \text{elsewhere}.
%%     \end{cases}
%% \end{equation}
%% In this case, $p=1$ and the coefficients $h_i$ and $g_i$ are $h_1=h_2=\frac{\sqrt{2}}{2}$ and
%% $g_1=\frac{\sqrt{2}}{2}$, $g_2=-\frac{\sqrt{2}}{2}$.
%% \end{example}
\begin{equation} \label{eq:cell_psi}
  \psi\wl{\ell}{k}(x)    = 2^{\ell/2}\psi\left( 2^{\ell}x-k \right),
  \quad
  \ell=0,1,2,\dots,\quad k\in\bbZ.
\end{equation}
For a given function $v(x)$ defined in $\mathbb{R}$, its scaling and wavelet coefficients are given by
\begin{equation}
  \label{eq:def_vd}
  v\wl{\ell}{k} := \int v(x)\varphi\wl{\ell}{k}(x)\dd x,
  \quad
  d\wl{\ell}{k} := \int v(x)\psi\wl{\ell}{k}(x)\dd x,
  \quad
  \ell=0,1,2,\dots,\quad k\in\bbZ.
\end{equation}
%% Note that the support of $\varphi\wl{\ell}{k}(x)$ maybe not in the domain $\Omega=[0,1)$. A
%%   technically convenient way is to pad the function $v(x)$ periodically. Hereafter, if not
%%   specified, we always assume the functions are periodic.
The refinement equation \eqref{eq:recursive_phi} implies that
\begin{equation}\label{eq:v2v}
  \begin{aligned}
    v\wl{\ell}{k} &= \int v(x) 2^{\ell/2}\varphi(2^{\ell}x-k) \dd x  =
    \int v(x) 2^{(\ell+1)/2}\sum_{i\in\bbZ} h_i\varphi(2^{\ell+1}x-2k-i) \dd x \\ 
    &= \int v(x) \sum_{i\in\bbZ} h_i\varphi\wl{\ell+1}{2k+i}(x)\dd x       = \sum_{i\in\bbZ} h_i v\wl{\ell+1}{2k+i},
  \end{aligned}
\end{equation}
which is a recursive relation between the coefficients $v\wl{\ell}{k}$ and $v\wl{\ell+1}{k}$. A
similar relation can be derived between $d\wl{\ell}{k}$ and $v\wl{\ell+1}{k}$:
\begin{equation}\label{eq:v2d}
  d\wl{\ell}{k}=\sum_{i\in\bbZ}g_i v\wl{\ell+1}{2k+i}.
\end{equation}
If $v(x)\in L^2(\bbR)$, then $v\sps{\ell}:=\left(v\wl{\ell}{k}\right)_{k\in\bbZ}$ and
$d\sps{\ell}:=\left(d\wl{\ell}{k}\right)_{k\in\bbZ}$ are sequences in $\ell^2(\bbZ)$. In the
operator form, \eqref{eq:v2v} and \eqref{eq:v2d} can be written as
\begin{equation}\label{eq:vd_matrix}
  v\sps{\ell}=(H\sps{\ell})^T v\sps{\ell+1},\quad 
  d\sps{\ell}=(G\sps{\ell})^T v\sps{\ell+1},
\end{equation}
where the operators $H\sps{\ell}, G\sps{\ell}: \ell^2(\bbZ)\rightarrow \ell^2(\bbZ)$ are banded with
width $2p$ due to the support size of $\{h_i\}$ and $\{g_i\}$, respectively. The orthogonality
relation \eqref{eq:orthogonality} can now be written compactly as
\begin{equation}
  (H\sps{\ell})^T H\sps{\ell} = I,\quad
  (G\sps{\ell})^T G\sps{\ell} = I,\quad
  H\sps{\ell} (H\sps{\ell})^T + G\sps{\ell} (G\sps{\ell})^T=I.
\end{equation}
By introducing an orthogonal operator $W\sps{\ell} = \begin{pmatrix} G\sps{\ell} & H\sps{\ell} \end{pmatrix}$,
one can rewrite \eqref{eq:vd_matrix} as
\begin{equation}\label{eq:d2v_v2d}
  \begin{pmatrix}
    d\sps{\ell}\\
    v\sps{\ell}
  \end{pmatrix}
  = (W\sps{\ell})^T v\sps{\ell+1}, \quad
  v\sps{\ell+1} = W\sps{\ell}
  \begin{pmatrix}
        d\sps{\ell}\\
        v\sps{\ell}
    \end{pmatrix},
    \quad \ell=0, 1,2,\dots,
\end{equation}
which are the decomposition and reconstruction steps of the wavelet analysis, respectively.  This
process can be illustrated by the following diagram
\begin{equation}
    \newcommand\lra{\longrightarrow}
    \setlength{\arraycolsep}{2pt}
    \begin{array}{ccccccccccccccc}
        \cdots & \lra & v\sps{\ell}& \lra     & v\sps{\ell-1} & \lra     & v\sps{\ell-2} & \lra
            & \cdots   & \lra     & v\sps{2} & \lra       & v\sps{1} & \lra      & v\sps{0} \\
        & \searrow &      & \searrow &          & \searrow &            & \searrow
            &          & \searrow &          & \searrow   &          & \searrow \\
        & &  d\sps{\ell}    &          & d\sps{\ell-1} &          & d\sps{\ell-2} &
            & \cdots   &          & d\sps{2} &            & d\sps{1} &           & d\sps{0}
    \end{array}.
\end{equation}

Although, until this point our discussion is concerned with wavelets defined on $\bbR$, the
definitions can be easily extended to the case of periodic functions defined on a finite interval,
say $[0,1]$. The only modification is that all the shifts and scalings in the $x$ variable are now
done modulus the integers. In addition, instead of going all the way to level $0$, the above
multiresolution analysis stops at a coarse level $L_0=O(\log_2 p)$ before the wavelet and scaling
functions start to overlap itself, which is diagramed as
\begin{equation}
    \newcommand\lra{\longrightarrow}
    \setlength{\arraycolsep}{2pt}
    \begin{array}{ccccccccccccc}
        \cdots & \lra &
      v\sps{\ell}& \lra     & v\sps{\ell-1} & \lra     & v\sps{\ell-2} & \lra
              & \cdots   & \lra     & v\sps{L_0} \\
        & \searrow &      & \searrow &          & \searrow &            & \searrow
              &          & \searrow &          \\
              & & d\sps{\ell}      &          & d\sps{\ell-1} &          & d\sps{\ell-2} &
              & \cdots   &          & d\sps{L_0} 
    \end{array}.
\end{equation}

Finally, the Daubechies wavelet $\psi(x)$ described above has exactly $p$ vanishing moments
\cite{daubechies1988orthonormal}, i.e.,
\begin{equation}
  \label{eq:vm}
  \int x^j \psi(x) \dd x=0, \quad j=0,\dots,p-1.
\end{equation}
Since the space of polynomials up to degree $p-1$ is invariant under scaling and translation,
\eqref{eq:vm} is true also for any wavelet $\psi^{(\ell)}_k(x)$.

%% To calculate the scaling and wavelet coefficients of a function, we proceed as follows. Suppose we
%% have $N=2^L$ samples $v\in\bbR^N$ of a function $v(x)$, which can for simplify be thought 
%% of as $v\sps{L}$ or obtained by some other discretization. Let $v\sps{L}=v$, then all the scaling
%% and wavelet coefficients can be obtained by evaluating
%% \begin{equation}
%%     \begin{pmatrix}
%%         d\sps{\ell}\\
%%         v\sps{\ell}
%%     \end{pmatrix}
%%     = W\sps{\ell} v\sps{\ell+1}, \quad
%%     \ell=L-1,L-2,\dots,0.
%% \end{equation}

%============
\subsection{Integral operator compression}\label{sec:wl2d}

The nonstandard form is concerned with wavelet based compression of integral operators. Suppose that $A$
is an integral operator on the periodic interval $[0,1]$ with kernel $a(x,y)$. We consider the
Galerkin projection of $A$ in the space $V_{L}$ for a sufficient deep level $L$. The $2^L\times 2^L$
matrix $A\sps{L} = (A\sps{L}_{k_1,k_2})_{k_1,k_2=0,\dots,2^L-1}$ with entries given by
\[
A\wl{L}{k_1,k_2} :=\iint \varphi\wl{L}{k_1}(x)a(x,y)\varphi\wl{L}{k_2}(y)\dd x\dd y
\]
is a reasonably accurate approximation of the operator $A$. The nonstandard form is essentially a
data-sparse representation of $A\sps{L}$ using the 2D multiresolution wavelet basis. Let us
now introduce
\begin{equation}\label{eq:Int_D2A}
  \begin{aligned}
      D\wl{\ell}{1,k_1,k_2}&:=\iint \psi\wl{\ell}{k_1}(x)a(x,y)\psi\wl{\ell}{k_2}(y)\dd x\dd y,
      &
      D\wl{\ell}{2,k_1,k_2}&:=\iint \psi\wl{\ell}{k_1}(x)a(x,y)\varphi\wl{\ell}{k_2}(y)\dd x\dd y,\\
      D\wl{\ell}{3,k_1,k_2}&:=\iint \varphi\wl{\ell}{k_1}(x)a(x,y)\psi\wl{\ell}{k_2}(y)\dd x\dd y,
      &
      A\wl{\ell}{k_1,k_2}  &:=\iint \varphi\wl{\ell}{k_1}(x)a(x,y)\varphi\wl{\ell}{k_2}(y)\dd x\dd y,
    \end{aligned}
\end{equation}
for $\ell=L_0, \dots, L-1$, and $k_1,k_2 = 0,\dots,2^\ell-1$. In what follows, it is convenient to
organize these coefficients into the following matrices:
\begin{equation}\label{eq:AD}
  A\sps{\ell}=(A\sps{\ell}_{k_1,k_2})_{k_1,k_2=0,\dots,2^\ell-1},\quad\quad
  D_j\sps{\ell}=(D\sps{\ell}_{j,k_1,k_2})_{k_1,k_2=0,\dots,2^\ell-1}, j=1,2,3.
\end{equation}
Using \eqref{eq:d2v_v2d}, these matrices can be computed from level $L$ down to level $L_0$ via the
recursive relation
\begin{equation}\label{eq:D2A}
  \begin{pmatrix}
    D_1\sps{\ell}   & D_2\sps{\ell}\\
    D_3\sps{\ell}   & A\sps{\ell}
  \end{pmatrix}
  = (W\sps{\ell})^T A\sps{\ell+1} W\sps{\ell}.
  %% A\sps{\ell+1} = W\sps{\ell}
  %% \begin{pmatrix}
  %%   D_1\sps{\ell}   & D_2\sps{\ell}\\
  %%   D_3\sps{\ell}   & A\sps{\ell}
  %% \end{pmatrix} (W\sps{\ell})^T,
  %% \quad \ell=0,\dots,L-1,
\end{equation}

%For a two-dimensional function $a(x,y)$, there
%are four part: wavelet-wavelet, wavelet-scaling, scaling-wavelet and scaling-scaling parts as

%\LY{is the following redundant?}
%Similar as one-dimensional case, to obtain all the scaling and wavelet coefficients of a
%two-dimensional function, we suppose that $N^2$ with $N=2^L$ samples $A\in\bbR^{N\times N}$ of a
%function $a(x,y)$ are given, which can for simplify be thought of as $A\sps{L}$ or obtained by some
%other discretization. Then all the scaling and wavelet coefficients can be obtained by letting
%$A\sps{L}=A$ and evaluating 
%\begin{equation}
%    \begin{pmatrix}
%        D_1\sps{\ell}   & D_2\sps{\ell}\\
%        D_3\sps{\ell}   & A\sps{\ell}
%    \end{pmatrix}
%    = (W\sps{\ell})^T A\sps{\ell+1} W\sps{\ell},\quad \quad \ell=L-1,L-2,\dots,0.
%\end{equation}

A key feature of \eqref{eq:D2A} is that the entries of the matrices $D_{j}\sps{\ell} (j=1,2,3)$
decay rapidly away from the diagonal if the kernel function $a(x,y)$ satisfies certain smoothness
conditions. More precisely, if $A$ is a Calderon-Zygmund operator with kernel $a(x,y)$ that
satisfies
\begin{equation}
  |a(x,y)|\lesssim \frac{1}{|x-y|},\quad 
  |\partial_x^{p} a(x,y)| + |\partial_y^{p} a(x,y)| \lesssim_p \frac{1}{|x-y|^{1+{p}}},
\end{equation}
then the entries of $D_j\sps{\ell}$ satisfy the following estimate \cite[Proposition 4.1]{bcr}
\begin{equation}
  |D_{j,k_1,k_2}\sps{\ell}| \lesssim_p \frac{C_{p}}{1+|k_1-k_2|^{p+1}},\quad \text{ for all } j
  \text{ and } |k_1-k_2|\geq 2p.
\end{equation}
Similarly, suppose that $A$ is a pseudo-differential operator with symbol $\sigma(x,\xi)$
\begin{equation}
  (Av)(x) = \int e^{\mathrm{i} x\xi}\sigma(x,\xi)\hat{v}(\xi)\dd\xi=\int a(x,y)v(y)\dd y.
\end{equation}
%where $a(x,y)$ is the distributional kernel of $A$, and $\hat{v}(\xi)$ is the Fourier transform of $v(x)$. 
If the symbols $\sigma$ of $A$ and $\sigma^*$ of $A^*$ satisfy the conditions
\begin{equation}
  |\partial_{\xi}^{m_1}\partial_{x}^{m_2}\sigma(x,\xi)|   \lesssim_{m_1,m_2}(1+|\xi|)^{\lambda-m_1+m_2},
  \quad
  |\partial_{\xi}^{m_1}\partial_{x}^{m_2}\sigma^*(x,\xi)| \lesssim_{m_1,m_2}(1+|\xi|)^{\lambda-m_1+m_2},
\end{equation}
then \cite[page 155]{bcr}
\begin{equation}
  \frac{1}{2^{\lambda \ell}} |D_{j,k_1,k_2}\sps{\ell}| \lesssim_p
  \frac{1}{(1+|k_1-k_2|)^{p+1}},\quad \text{ for all } j \text{ and } k_1, k_2.
\end{equation}
Because of the rapid decay of the entries of $D_{j,k_1,k_2}\sps{\ell}$, one can approximate these
matrices by truncating at a band of width $n_b = O(\log(1/\eps))$ for a prescribed {\em relative
  accuracy} $\eps$.  Note that, the value of $n_b$ is independent of the specific choices of
$\ell=L_0,\dots,L-1$, $j=1,2,3$, or the mesh size $N$. From now on, we assume that the matrices
$D_{j,k_1,k_2}\sps{\ell}$ are pre-truncated already. The number of non-zero entries of these
matrices at level $\ell$ is clearly $O(2^\ell)$.

%In the following, we often approximate the matrices $D\sps{\ell}_j$ by a cyclic band matrices
%\fy{?} and we abuse the notation calling the band matrices $D\sps{\ell}_j$. The band size of
%$D_{j}\sps{\ell}$ for all $\ell=0,\dots,L$ and $j=1,2,3$ is denoted by $n_b$, which is a constant
%independent of mesh size $N$.

Concerning the matrix $A\sps{\ell}$ for $\ell=L_0, \dots, L-1$, though they are dense in general,
the recursive relation \eqref{eq:AD} shows that one only needs to store the matrix $A\sps{L_0}$ at
the top level $L_0$. Since there are only $O(p)$ wavelets and scaling functions at level $L_0$,
$A\sps{L_0}$ is of constant size. Combining this with the estimate for matrices
$D_{j,k_1,k_2}\sps{\ell}$ over all levels, one concludes that the total number of non-zero
coefficients in the nonstandard form is $O(N)$. This is the {\em nonstandard form} proposed by
Beylkin, Coifman, and Rokhlin.

%============
\subsection{Matrix-vector multiplication}

For an integral operator $A$, a fundamental operation is to apply $A$ to an arbitrary function $v$,
i.e.,
\begin{equation}\label{eq:integral}
  u = A v, \quad  u(x) = \int_{\Omega}a(x,y)v(y) \dd y, \quad \Omega=[0,1],
\end{equation}
where $v$ and $u$ are periodic functions defined on $\Omega$. With the Galerkin discretization
described above at level $L$, this is simply a matrix vector multiplication.
\begin{equation}\label{eq:discrete}
  u\sps{L} = A\sps{L} v\sps{L}.
\end{equation}

\begin{figure}[ht]
  \centering
  \subfloat[wavelet transform]{
    \includegraphics[height=0.18\textheight,clip]{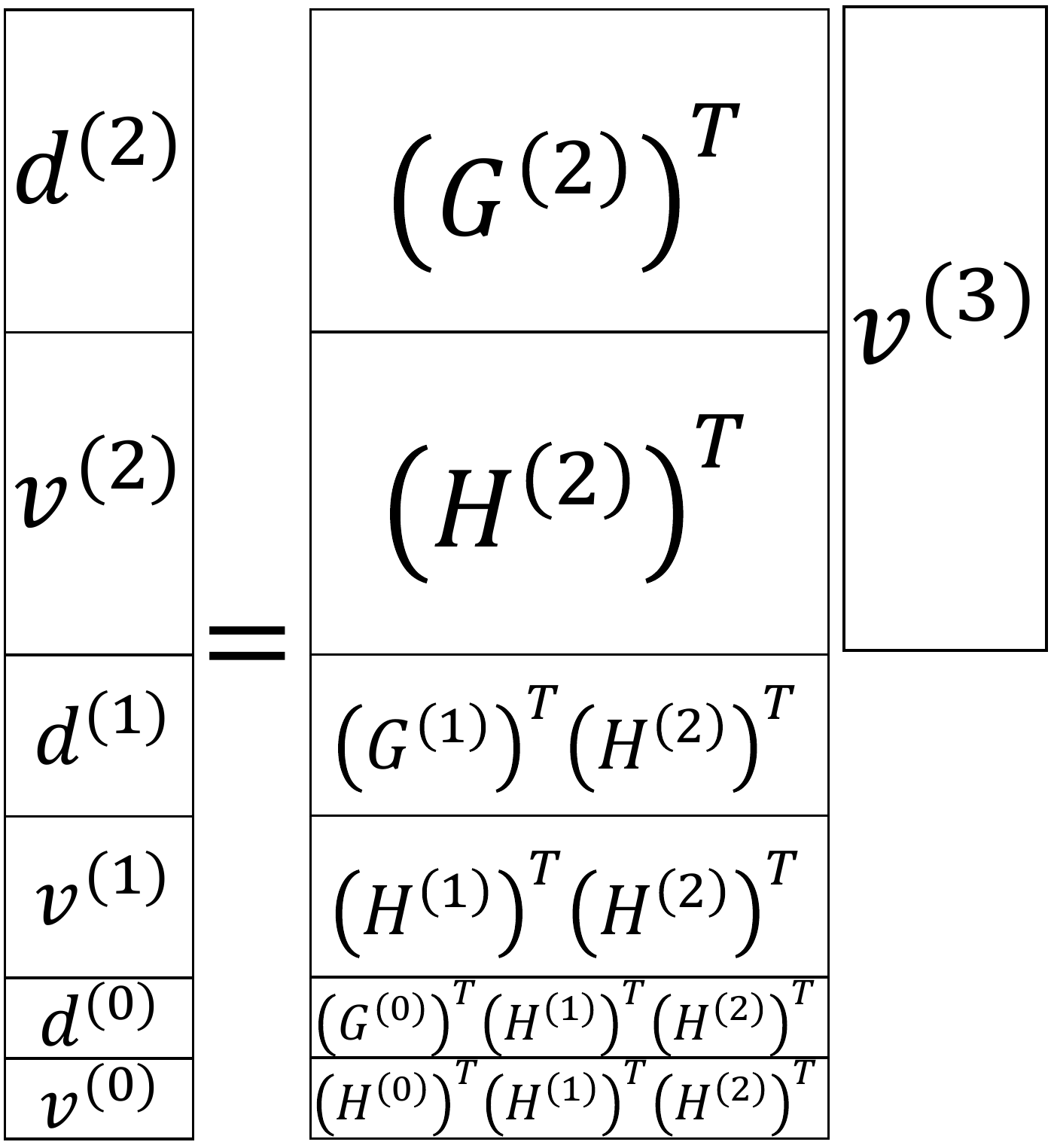}
  }~
  \subfloat[matrix-vector multiplication]{
    \includegraphics[height=0.18\textheight,clip]{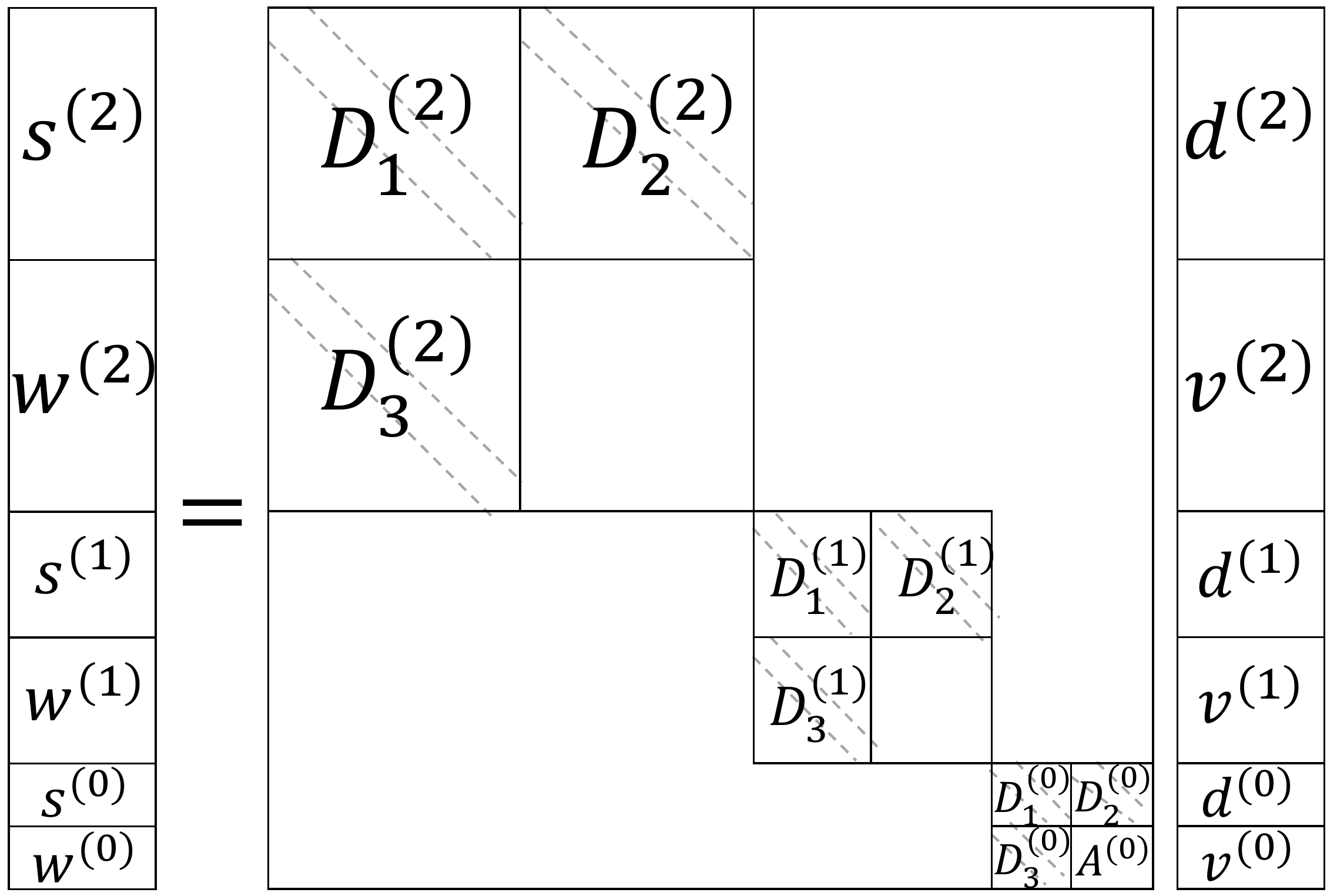}
  }~
  \subfloat[inverse wavelet transform]{
    \includegraphics[height=0.18\textheight,clip]{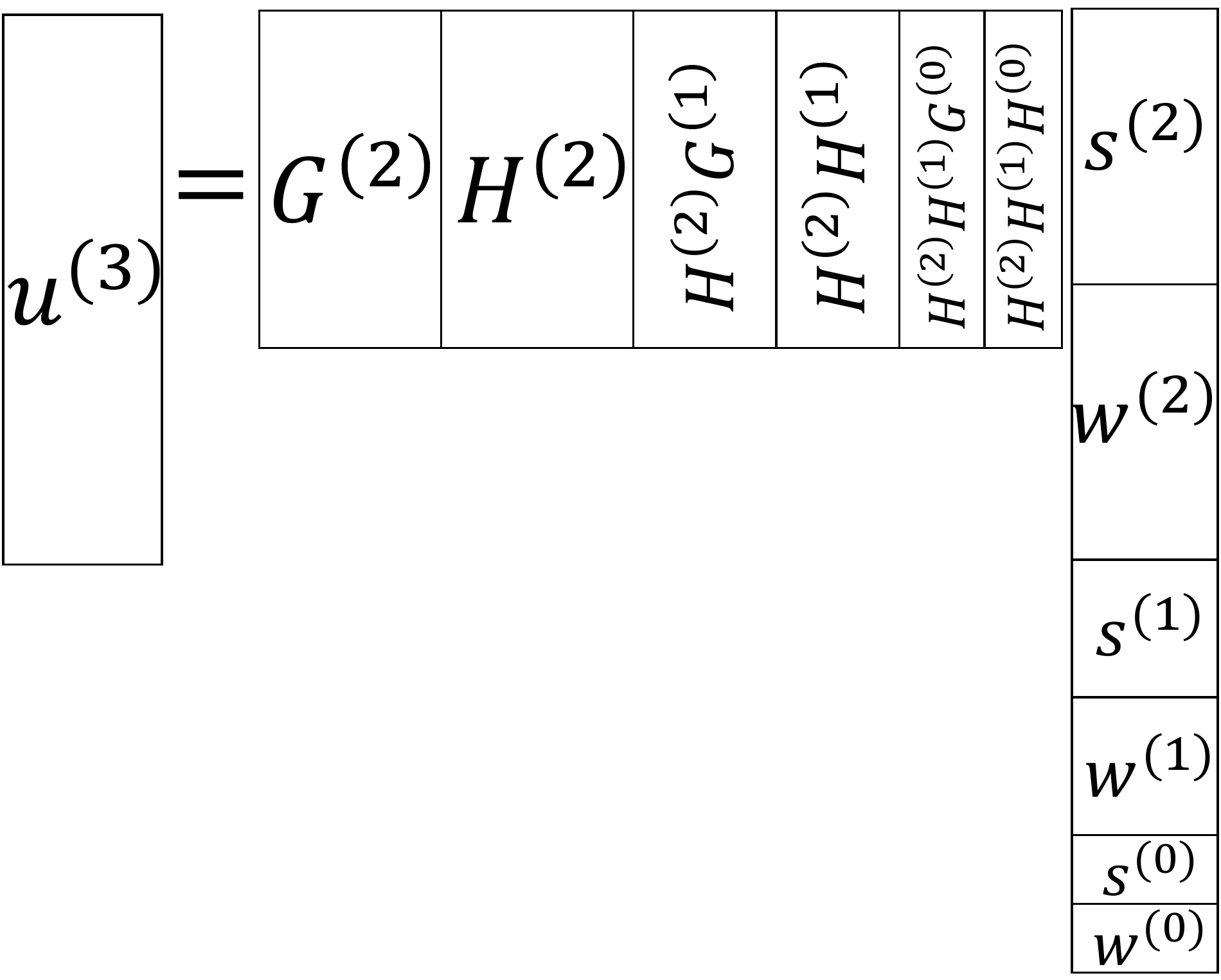}
  }
  \caption{\label{fig:matvec}Illustration of matrix-vector multiplication based on the nonstandard
    form (with $L_0=0$ and $L=3$). $D\sps{\ell}_j$, $j=1,2,3$ are all band matrices.}
\end{figure}

This matrix-vector multiplication is computationally intensive since $A\sps{\ell}$ is a dense
matrix. However, the nonstandard form introduced above offers an $O(N)$ (linear complexity)
algorithm for carrying out the matrix vector multiplication. The key observation of this process is
the following recurrence relation
\begin{equation}\label{eq:A2D}
  A\sps{\ell+1} = W\sps{\ell}
  \begin{pmatrix}
    D_1\sps{\ell}   & D_2\sps{\ell}\\
    D_3\sps{\ell}   & A\sps{\ell}
  \end{pmatrix} (W\sps{\ell})^T,
  \quad \ell=0,\dots,L-1,
\end{equation}
obtained from transposing \eqref{eq:D2A}. This means that, modulus wavelet transforms, at each level
$\ell$, one only needs to perform matrix-vector multiplication with the sparse matrices
$D_{j,k_1,k_2}\sps{\ell}$. The dense matrix multiplication is only carried out with $A\sps{L_0}$ at
the top level $L_0$. More precisely, the algorithm follows the three steps illustrated in
\cref{fig:matvec}:
\begin{enumerate}
\item[(a)] Transforming $v\sps{L}$ to obtain its scaling and wavelet coefficients of the
  nonstandard form;
\item[(b)] Multiplying with banded matrices at each scale;
\item[(c)] Transforming the scaling and wavelet coefficients of the nonstandard form to obtain
  $u\sps{L}$.
\end{enumerate}

\begin{figure}[ht]
    \centering
    \includegraphics[page=1,width=\textwidth]{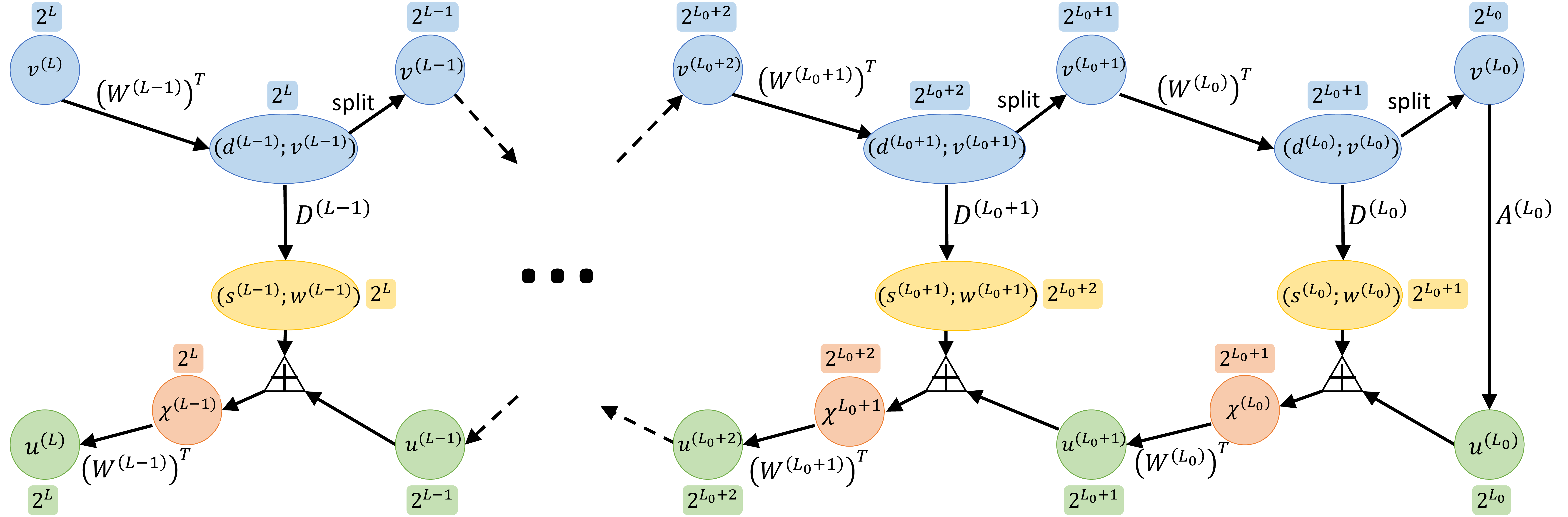}
    \newsavebox{\wlw}
    \savebox{\wlw}{$\chi\sps{\ell}=\begin{pmatrix}
        s\sps{\ell}\\
        w\sps{\ell}
    \end{pmatrix}+\begin{pmatrix}
        0\\
        u\sps{\ell}
    \end{pmatrix}$}
    \caption{\label{fig:wl}Diagram of the matrix-vector multiplication in the nonstandard form. The
      notation $\trianglesignplus$ denotes the summation \usebox{\wlw}.}
\end{figure}

By following \eqref{eq:d2v_v2d}, \eqref{eq:D2A}, \eqref{eq:A2D} and introducing auxiliary vectors
$u\sps{\ell}$ at all levels, this calculation can be written more compactly as
\begin{equation}\label{eq:matvec}
  u\sps{\ell+1} = A\sps{\ell+1}v\sps{\ell+1}
  %= (W\sps{\ell})^T \begin{pmatrix}
  %    D_1\sps{\ell} & D_2\sps{\ell}\\
  %    D_3\sps{\ell} & A\sps{\ell}
  %\end{pmatrix} W\sps{\ell}
  %(W\sps{\ell})^T \begin{pmatrix}
  %    d\sps{\ell}\\
  %    v\sps{\ell}
  %\end{pmatrix}
  = W\sps{\ell} \begin{pmatrix}
    D_1\sps{\ell} & D_2\sps{\ell}\\
    D_3\sps{\ell} & A\sps{\ell}
  \end{pmatrix}
  \begin{pmatrix}
    d\sps{\ell}\\
    v\sps{\ell}
  \end{pmatrix}
  = W\sps{\ell} \left[\begin{pmatrix}
      D_1\sps{\ell} & D_2\sps{\ell}\\
      D_3\sps{\ell} & D_4\sps{\ell}
    \end{pmatrix}
    \begin{pmatrix}
      d\sps{\ell}\\
      v\sps{\ell}
    \end{pmatrix}
    + \begin{pmatrix}
      0\\
      u\sps{\ell}
    \end{pmatrix}\right],
\end{equation}
for $\ell=L-1,L-2,\dots,L_0$ with $D_4\sps{\ell}$ defined to be the zero matrix. \cref{alg:wl}
provides the pseudo code for computing \eqref{eq:matvec}. \cref{fig:wl} illustrates the
multiresolution structure of the computation, which will turn out to be useful when we introduce the
NN. Since the matrices $D_j\sps{\ell}$ are band matrices, the arithmetic complexity of evaluating
$u\sps{\ell+1}$ given $u\sps{\ell}$ is $O(n_bN / 2^\ell)$, where $N=2^L$. Therefore, the overall
complexity of \eqref{eq:discrete} is $O(N)$.

\begin{algorithm}[htb]
\begin{small}
\begin{center}
\begin{minipage}{0.5\textwidth}
\begin{algorithmic}[1]
  %\State $v\sps{L}=v$;
  \For {$\ell$ from $L-1$ to $L_0$ by $-1$}\vspace{2mm}
  \State $\begin{pmatrix}
    d\sps{\ell}\\
    v\sps{\ell}
  \end{pmatrix}=(W\sps{\ell})^T v\sps{\ell+1}$;\vspace{2mm}
  \EndFor
  \State $u\sps{L_0}=A\sps{L_0}v\sps{L_0}$;
  \algstore{break1}
\end{algorithmic}
\vspace{-4mm}
\end{minipage}
\begin{minipage}{0.4\textwidth}
\begin{algorithmic}[1]
  \algrestore{break1}
  \For {$\ell$ from $L_0$ to $L-1$}
  \State {\footnotesize$
    \begin{pmatrix}
        s\sps{\ell}\\
        w\sps{\ell}
    \end{pmatrix}
    =\begin{pmatrix}
    D_1\sps{\ell} & D_2\sps{\ell} \\
    D_3\sps{\ell} & D_4\sps{\ell}
    \end{pmatrix}
    \begin{pmatrix}
      d\sps{\ell}\\
      v\sps{\ell}
    \end{pmatrix}$};
  \State {\footnotesize $u\sps{\ell+1}= W\sps{\ell} \begin{pmatrix}
      s\sps{\ell}\\
      w\sps{\ell}+u\sps{\ell}
    \end{pmatrix}$};
  \EndFor
\end{algorithmic}
\vspace{-4mm}
\end{minipage}
\end{center}
\end{small}
\caption{$u\sps{L} = A\sps{L} v\sps{L}$ in the nonstandard form.}
\label{alg:wl}
\end{algorithm}

%============
\subsection{Multidimensional case}

The above discussion can be easily extended to higher dimensions by using the multidimensional
orthonormal wavelets. The discussion here will focus on the two-dimensional case. For the
two-dimensional wavelet analysis (see for example \cite[Chapter 8]{misiti2013wavelets} and
\cite{Mallat2009}), there are three different types of wavelets at each scale.

The 2D analog of the recursive relation \eqref{eq:A2D} takes the form
\begin{equation}\label{eq:D2A2D}
  \begin{pmatrix}
    D_1\sps{\ell}   & D_2\sps{\ell}   & D_3\sps{\ell} & D_4\sps{\ell}\\
    D_5\sps{\ell}   & D_6\sps{\ell}   & D_7\sps{\ell} & D_8\sps{\ell}\\
    D_9\sps{\ell}   & D_{10}\sps{\ell} & D_{11}\sps{\ell} & D_{12}\sps{\ell}\\
    D_{13}\sps{\ell} & D_{14}\sps{\ell}  & D_{15}\sps{\ell} & A\sps{\ell}
  \end{pmatrix}
  = (W\sps{\ell})^T A\sps{\ell+1} W\sps{\ell}
\end{equation}
and all $D\sps{\ell}_j$ matrices are sparse with only $O(2^\ell)$ non-negligible entries after
thresholding. The matrix-vector multiplication takes exactly the same form as \cref{alg:wl} except
obvious changes due to matrix sizes.

%channels in the neural network changes from 2 in 1D to $2^2$ in 2D. We adopt the convention that the
%first three channels correspond to the wavelet coefficients and the last one to the scaling function
%coefficients. Therefore,

%=======================================================
\section{Matrix-vector multiplication as a neural network}\label{sec:mvnn}

The goal of this section is to represent the matrix-vector multiplication in \cref{alg:wl} as a
linear neural network. We start by introducing several basic tools in \cref{sec:cnn} and then
present the neural network representation of \cref{alg:wl} in \cref{sec:NNwl}. The presentation will
be focused on the 1D case first in order to illustrate the main ideas clearly.

%==========
\subsection{Neural network tools}\label{sec:cnn}

Throughout the discussion below, the input, output, and intermediate data are all represented with
2-tensor. For a 2-tensor of size $N_x\times \alpha$, we call $N_x$ the \emph{spatial dimension} and
$\alpha$ the \emph{channel dimension}.

To perform the operations appeared in \cref{alg:wl}, we first introduce a few common NN layers. The
first one is the well-known \emph{convolutional layer} (Conv) where the output of each location
depends only locally on the input. Given an input tensor $\xi$ of size $N_{\ipt}\times\alpha_{\ipt}$
and an output tensor $\zeta$ of size $N_{\out}\times\alpha_{\out}$, the convolutional layer performs
the computation
\begin{equation}\label{eq:conv}
  \zeta_{i,c'} = \phi\left( \sum_{j=is}^{is+w-1}\sum_{c=0}^{\alpha_{\ipt}-1}
  W_{j;c',c}\xi_{j,c} + b_{c'}\right),
  \quad i=0,\dots,N_{\out}-1, ~ c'=0,\dots,\alpha_{\out}-1,
\end{equation}
where $\phi$ is a pre-specified function, called \emph{activation}, usually chosen to be a linear
function, a rectified-linear unit (ReLU) function, or a sigmoid function. The parameters $w$ and $s$
are called the \emph{kernel window size} and \emph{stride}, respectively. Here we assume that
$N_{\out}=N_{\ipt}/s$ and the tensor $\xi$ is periodically padded if the index is out of range. Note
that this differs from the definition of the convolution layer in TensorFlow \cite{tensorflow} for
which zero padding is the default behavior. In what follows, such a convolution layer is denoted as
\begin{equation}
  \zeta = \Conv[\alpha_{\out}, w, s, \phi](\xi),
\end{equation}
where the values of $N_{\ipt}$, $\alpha_{\ipt}$, and $N_{\out}$ are inferred from the input
tensor $\xi$.

%We call the layer $\zeta$ \emph{convolutional layer} (\Conv~layer) hereafter.
%\LY{comment that here we assume that things are periodic and also the convolution window is slightly
%different from Keras/Tensorflow. But it does not cause any trouble. Just a shift}

Note that the weight $W_{j;c',c}$ in \eqref{eq:conv} is independent on the position $i$ of $\zeta$,
thus \Conv is translation invariant. When the weight are required to depend on the position $i$
(i.e. $W_{i,j;c',c}$), the natural extension of \Conv is the so-called \emph{locally connected} (LC)
layer. This layer is denoted by
\[
\zeta = \LC[\alpha_{\out},w,s,\phi](\xi).
\]

Finally, when the output data tensor depends on every entry of the input tensor, this is called a
{\em dense} layer, denoted by
\[
\zeta = \Dense[\alpha_{\out},\phi](\xi).
\]
Here we assume implicitly that the spatial dimensions of the input and output tensors are same,
i.e., $N_{\out}=N_{\ipt}$.

% In \eqref{eq:conv} the convolutional network is represented using tensor notation; however, we can
% reshape $\zeta$ and $\xi$ to a vector by column major indexing and $W$ to a matrix and write
% \eqref{eq:conv} into a matrix-vector form as
% \begin{equation}\label{eq:convmat}
%     \zeta = \phi(W\xi+b).
% \end{equation}

%% The second tool is the period padding function. For a tensor $\xi$ of size $N\times \alpha$,
%% we define the period padding function as 
%% \begin{equation}\label{eq:def_pad}
%%     \PeriodPad[n](\xi) = \begin{pmatrix}
%%         \xi((N-n):N, :); & \xi; & \xi(1:n, :)
%%     \end{pmatrix}.
%% \end{equation}
%% Here the MATLAB notation is used. \note{Matlab or python}{?}

% We also need a reshape function $\Reshape[w]$ which reshapes a 2-tensor of size $n_1\times n_2$ to
% a 2-tensor of size $wn_1 \times n_2/w$ by row major indexing. Here, we implicitly regard a vector
% of size $n$ as a 2-tensor of size $n\times 1$. 

With these basic tools, we can implement the key steps of \cref{alg:wl} in the NN framework:
\begin{itemize}
\item {Multiply $(W\sps{\ell})^T$ with a vector}: this step takes the form
  \begin{equation}\label{eq:def_WT}
    \zeta = \Conv[2, 2p, 2, \id](\xi),
  \end{equation}
  where $\id$ is the identity operator. The size of $\zeta$ is $M / 2 \times 2$ if the size of $\xi$
  is $M\times 1$. The convention adopted is that the first channel is for the wavelet coefficients
  and the second channel is for the scaling function coefficients.
\item {Multiply $\begin{pmatrix}
    D_1\sps{\ell} & D_2\sps{\ell} \\
    D_3\sps{\ell} & D_4\sps{\ell}
  \end{pmatrix}$ with a vector}: this step takes the form
  \begin{equation}\label{eq:def_LCK}
      \zeta = \LC[2, n_{b}, 1, \id](\xi).
  \end{equation}
  The size of $\zeta$ is $M\times2$ if the size of $\xi$ is $M\times2$. Notice that the width of LC
  layer corresponds to the band width of the banded matrices.
  %Analogously, we define $\Conv$ by replacing $\LC$ in \eqref{eq:def_LCK} by $\Conv$.
\item {Multiply $A\sps{L_0}$ with a vector}: this step takes the form
  \[
  \zeta = \Dense[1,\id](\xi).
  \]
  Both $\zeta$ and $\xi$ are of size $2^{Lc}\times 1$.
\item {Multiply $W\sps{\ell}$ with a vector}: this step first computes
  \begin{equation}\label{eq:def_IWT}
    \zeta = \Conv[2, p, 1, \id](\xi),
  \end{equation}
  followed by a reshape that goes through the channel dimension first. The output is of size
  $2M\times1$ if the input is $M\times2$.
\end{itemize}
In all these NN operations, the bias $b$ is set to be zero.

%==========
\subsection{Neural network representation}\label{sec:NNwl}

\begin{figure}[!h]
  \centering
    \includegraphics[page=2,width=\textwidth]{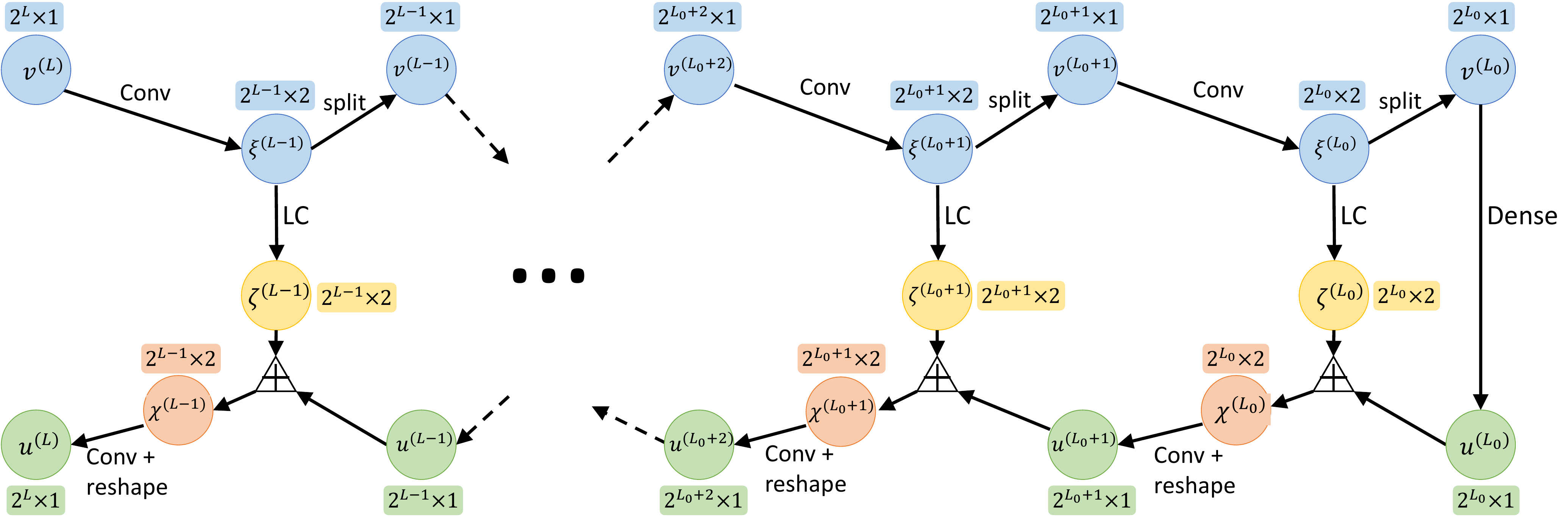}
    \newsavebox{\wlwn}
    \savebox{\wlwn}{$\xi\sps{\ell}+\begin{pmatrix}
      0\\
      u\sps{\ell}
      \end{pmatrix}$}
    \caption{\label{fig:nnwl} Neural network architecture of the matrix-vector multiplication in the
      nonstandard form.  ``split'' means extracting the last channel of $\xi\sps{\ell}$ to obtain
      $v\sps{\ell}$.  $\trianglesignplus$ means adding $u\sps{\ell}$ to the last channel of
      $\zeta\sps{\ell}$ to obtain $\chi\sps{\ell}$.}
\end{figure}

Combining the basic tools introduced above, we can translate \cref{alg:wl} into an NN. The resulting
algorithm is summarized in \cref{alg:nnwl}. \cref{fig:nnwl} illustrates the multiresolution
structure of the NN. It has the same structure as \cref{fig:wl} except the basic operations are
replaced with the NN layers introduced in \cref{sec:cnn}.
 
\begin{algorithm}[!h]
  \begin{small}
    \begin{center}
      \begin{minipage}{0.9\textwidth}
        \begin{algorithmic}[1]
          \For {$\ell$ from $L-1$ to $L_0$ by $-1$}
          \State $\xi\sps{\ell} = \Conv[2,2p,2,\id](v\sps{\ell+1})$; \label{alg:line_wt}
          \State $v\sps{\ell}$ is the last channel of $\xi\sps{\ell}$;
          \EndFor
          \State $u\sps{L_0}=\Dense[1,\id](v\sps{L_0})$;
          \For {$\ell$ from $L_0$ to $L-1$}
          \State $\zeta\sps{\ell}=\LC[2, n_b, 1,\id](\xi\sps{\ell})$;
          \State Adding $u\sps{\ell}$ to the last channel of $\zeta\sps{\ell}$ gives $\chi\sps{\ell}$;\label{alg:line_add}
          \State $u\sps{\ell+1}=\Conv[2, p, 1, \id](\chi\sps{\ell})$; \label{alg:line_iwt}
          \State Reshape $u\sps{\ell+1}$ to a 2-tensor of size $2^{\ell+1}\times1$ by following channel dimension first;
          \EndFor
        \end{algorithmic}
        \vspace{-4mm}
      \end{minipage}
    \end{center}
  \end{small}
  \caption{NN architecture for $u\sps{L} = A\sps{L} v\sps{L}$ in the nonstandard form.}
  \label{alg:nnwl}%
\end{algorithm}

Let us now count the number of parameters used in the network in \cref{alg:nnwl}. Notice that the
number of parameters in $\Conv$ and $\LC$ are $w\alpha_{\ipt}\alpha_{\out}$ and
$N_{\out}w\alpha_{\ipt}\alpha_{\out}$, respectively.  The total number of parameters in
\cref{alg:nnwl} is
\begin{equation}
  \Nparams = \sum_{\ell=L_0}^{L-1}\left( 4p + 2^{\ell}4n_b + 4p\right) + 4^{L_0} \approx 4Nn_b +
  8p(L-L_0) + 4^{L_0}.% = O(N).
\end{equation}
Since $2^{L_0}$ is a small constant, the total number of parameters is $O(Nn_b)$.

%==========
\subsection{Multidimensional case}\label{sec:nD}

Our discussion here focuses on the 2D case. Each piece of data in the 2D algorithm can be
represented by a 3-tensor of size $N_{x,1}\times N_{x,2}\times\alpha$, where $N_x=(N_{x,1},
N_{x,2})$ is the size in the spatial dimension and $\alpha$ is the channel number. If a tensor $\xi$
of size $N_{\ipt,1}\times N_{\ipt,2}\times \alpha_{\ipt}$ is connected to a tensor $\zeta$ of size
$N_{\out,1}\times N_{\out,2}\times \alpha_{\out}$ by a convolution layer, the computation takes the
form
\begin{footnotesize}
  \begin{equation}\label{eq:conv2d}
    \zeta_{i,c'} = \phi\left( \sum_{j_1=i_1s}^{i_1s+w-1}\sum_{j_2=i_2s}^{i_2s+w-1}\sum_{c=0}^{\alpha_{\ipt}-1}
    W_{j;c',c}\xi_{j,c} + b_{c'}\right),
    \quad i_1=0,\dots,N_{\out,1}-1, i_2=0, \dots,N_{\out,2}-1, c'=0,\dots,\alpha_{\out}-1.
  \end{equation}
\end{footnotesize}%
By denoting such a layer by $\Conv2$, we can write \eqref{eq:conv2d} concisely as
$\Conv2[\alpha,w,s,\phi](\xi)$. Similar to the 1D case, we also define the locally connected layer,
denoted by $\LC2[\alpha,w,s,\phi](\xi)$ and the dense layer, denoted by $\Dense2[\alpha,\phi](\xi)$.

With these tools, one can readily extend the \cref{alg:nnwl} to the 2D case, shown in
\cref{alg:nnwl2d}.
\begin{algorithm}[htb]
  \begin{small}
    \begin{center}
      \begin{minipage}{0.9\textwidth}
        \begin{algorithmic}[1]
          %\State $v\sps{L}=v$;
          \For {$\ell$ from $L-1$ to $L_0$ by $-1$}
          \State $\xi\sps{\ell} = \Conv2[4, 2p, 2,\id](v\sps{\ell+1})$; \label{alg:line_wtnn}
          \State $v\sps{\ell}$ is the last channel of $\xi\sps{\ell}$;
          \EndFor
          \State $u\sps{L_0}=\Dense2[1,\id](v\sps{L_0})$;
          \For {$\ell$ from $L_0$ to $L-1$}
          \State $\zeta\sps{\ell}=\LC2[4, n_b, 1, \id](\xi\sps{\ell})$;
          \State Adding $u\sps{\ell}$ to the last channel of $\zeta\sps{\ell}$ gives $\chi\sps{\ell}$;
          \State $u\sps{\ell+1}=\Conv2[4, p, 1, \id](\chi\sps{\ell})$; \label{alg:line_iwtnn}
          \State Reshape $u\sps{\ell+1}$ to a $3$-tensor of size $2\sps{\ell+1}\times
          2\sps{\ell+1}\times 1$;
          \EndFor
        \end{algorithmic}
        \vspace{-4mm}
      \end{minipage}
    \end{center}
  \end{small}
  \caption{NN architecture for $u\sps{L} = A\sps{L} v\sps{L}$ in the nonstandard form for the 2D
    case.}
  \label{alg:nnwl2d}%
\end{algorithm}

The reshape step at the end of \cref{alg:nnwl2d} deserves some comments. The input tensor
$u\sps{\ell+1}$ is of size $2^\ell\times 2^\ell \times 4$. The reshape process first change the
input to a $2^\ell\times 2^\ell \times 2 \times 2$ tensor by splitting the last dimension. It then
permutes the second and their third dimension to obtain a 4-tensor of size $2^\ell\times 2 \times
2^\ell \times 2$. By grouping the first and the second dimensions as well as the third and the
fourth dimension, one obtains a $2^{\ell+1}\times 2^{\ell+1}$ tensor. Finally, this tensor is
regarded as a 3-tensor of size $2^{\ell+1}\times 2^{\ell+1}\times 1$, i.e., with a single component
in the channel (last) dimension.

%% 3-tensor $\xi$ of size $n_1\times
%% n_2\times w^2\alpha$ to a 3-tensor $\zeta$ of size $n_1w\times n_2w\times \alpha$ by reshaping
%% $\xi_{j,k,:}$ to a 3-tensor of size $w\times w\times \alpha$ and joining them to a large 3-tensor.
%% %which is diagramed in \cref{fig:reshape2d}.

%==================================================================================================
\section{BCR-Net}\label{sec:bcrnet}

%%  In a standard approach for such
%% nonlinear maps, one uses iterative methods, which may require a large number of iterations and, at
%% each iteration, solve the linearized equation several times, resulting in computationally expensive
%% algorithms.

For many nonlinear map of form
\begin{equation}\label{eq:nonlinearmap}
    u=\cM(v),\quad  u,v\in\bbR^{N^d},
\end{equation}
when the singularity of $u$ only appears at the singularity of $v$, such a operator can be viewed as
a nonlinear generalization of pseudo-differential operators.  By leveraging the representation power
of the neural networks, we extend the architecture constructed in \cref{alg:nnwl} to represent the
nonlinear maps \eqref{eq:nonlinearmap}. The resulting NN architecture is referred to as the BCR-Net.
In order to simplify the presentation, we will focus on the 1D case in this section as all the results
here can be easily extended to the multi-dimensional case by following the discussion in
\cref{sec:nD}.

%\subsection{Algorithm and architecture}\label{sec:mnnAlgorithm}
Two changes are made in order to extend \cref{alg:nnwl} to the nonlinear case. The first is to
replace some of the identity activation functions $\id$ with the nonlinear activation
functions. More precisely, the activation function in $\LC$ and $\Dense$ is replaced with a
nonlinear one (either being ReLU or Sigmoid function), denoted by $\phi$.  On the other hand, since
each \Conv layer corresponds to a step of the wavelet transform, its activation function is kept
as $\id$.

\begin{figure}[h!]
    \centering
    \includegraphics[width=\textwidth]{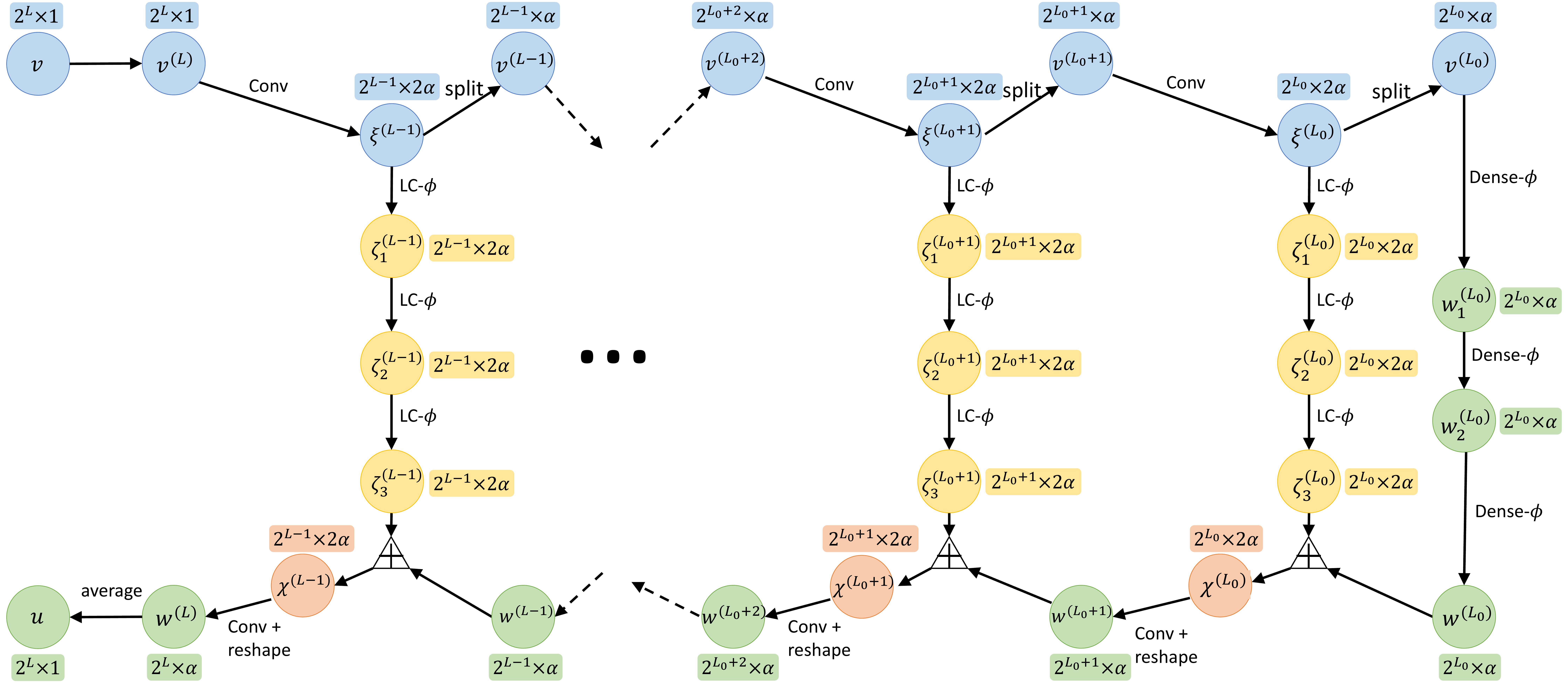}
    \caption{\label{fig:bcrnet} Architecture of BCR-Net.  ``split'' means extracting the last
      $\alpha$ channels of $\xi\sps{\ell}$ to obtain $v\sps{\ell}$.  $\trianglesignplus$ means
      adding $u\sps{\ell}$ to the last $\alpha$ channels of $\zeta\sps{\ell}$ to obtain
      $\chi\sps{\ell}$.}
\end{figure}

The second modification is to increase the ``width'' and ``depth'' of the network. The NN in
\cref{alg:nnwl} is rather narrow (the number of channels in $\LC$ is only $2$) and shallow (the
number of $\LC$ layers is only $1$). As a result, its representation power is limited. In order to
represent more general nonlinear maps, we increase the number of channels from $2$ to $2\alpha$,
with $\alpha$ wavelet channels and $\alpha$ scaling function channels. Here $\alpha$ is a
user-specified parameter. Moreover, the network also becomes deeper by increasing number of $\LC$
and $\Dense$ layers. The resulting algorithm is summarized in \cref{alg:bcrnet} and illustrated in
\cref{fig:bcrnet}.

\begin{algorithm}[htb]
\begin{small}
\begin{center}
  \newcommand\rankone{\alpha}
  \newcommand\ranktwo{\alpha}
  \begin{minipage}{0.9\textwidth}
    \begin{algorithmic}[1]
      \State $v\sps{L}=v$;
      \For {$\ell$ from $L-1$ to $L_0$ by $-1$}
      \State $\xi\sps{\ell} = \Conv[2\rankone,2p,2,\id](v\sps{\ell+1})$;
      \State $v\sps{\ell}$ is the last $\alpha$ channels of $\xi\sps{\ell}$;
      \EndFor
      \State $u_0\sps{L_0}=v\sps{L_0}$;
      \For {$k$ from $1$ to $K$ do}
      \State $u_{k}\sps{L_0}=\Dense[\ranktwo,\phi](u\sps{L_0}_{k-1})$;
      \EndFor
      \State $u\sps{L_0}=u_K\sps{L_0}$;
%      \algstore{break3}
%    \end{algorithmic}
%  \end{minipage}
%  \begin{minipage}{0.5\textwidth}
%    \begin{algorithmic}[1]
%      \algrestore{break3}
      \For {$\ell$ from $L_0$ to $L-1$}
      \State $\zeta_0\sps{\ell}=\xi\sps{\ell}$;
      \For {$k$ from $1$ to $K$}
      \State $\zeta\sps{\ell}_k=\LC[2\ranktwo, n_b, 1, \phi](\zeta\sps{\ell}_{k-1})$;
      \EndFor
      %\State $t_1 = \zeta\sps{\ell}_K(1:\ranktwo,:)$;
      %$\quad t_2 = \zeta\sps{\ell}_K(\ranktwo+1:2\ranktwo,:)$;
      %\State $t_3 = t_2 + w\sps{\ell}_K$;
      %\State $\chi\sps{\ell} = (\wc(\zeta\sps{\ell}_K),      \scc(\zeta\sps{\ell}_K)+u\sps{\ell})$;
      \State Adding $u\sps{\ell}$ to the last $\alpha$ channels of $\zeta\sps{\ell}_K$ gives $\chi\sps{\ell}$;\label{alg:line_add_mnn}
      \State $u\sps{\ell+1}=\Conv[2\rankone, p, 1, \id](\chi\sps{\ell})$; 
      \State Reshape $u\sps{\ell+1}$ to a 2-tensor of size $2^{\ell+1}\times\alpha$ by following channel dimension first;
      \EndFor
      \State Average over the channel direction of $u\sps{L}$ to give $u$;
    \end{algorithmic}
        \vspace{-4mm}
  \end{minipage}
\end{center}
\end{small}
\caption{BCR-Net applied to an input $v\in\bbR^N$ with $N=2^L$.}
\label{alg:bcrnet}%
\end{algorithm}

%% We remark that the function $\scc(\xi)$ and $\wc(\xi)$ in \cref{alg:bcrnet} extract the scaling and
%% wavelet parts, defined by
%% \begin{equation}
%%     \wc(\xi) = \xi(:,1:\alpha),\quad \scc(\xi)=\xi(:,\alpha+1:2\alpha).
%% \end{equation}
% To obtain a single value per location in $u$, we use the function $\ave(\xi)$ in \cref{alg:bcrnet} to calculate the weighted average of $\xi$ with respect to different channels.

% defined by
% \begin{equation}
%     \ave(\xi) = \sum_{i=1}^\alpha w_i\xi(:,i),
% \end{equation}
% where $\alpha$ is the number of channels of $\xi$ and $w_i$ are parameters to be determined.

%In the above, we focused on the case $N=2^L=2^{L-L_0}\times 2^{L_0}$. It is straightforward to
%extend it to the case $N=2^{L-L_0}m$ for any positive integer $m$.

Similar to the linear case, the number of parameters of BCR-Net can be estimated as follows:
\begin{equation}
  \begin{aligned}
    \Nparams &= 
    \sum_{\ell=L_0}^{L-1}\left(4\alpha^2 p 
    + \sum_{k=1}^K(2\alpha)^2n_b2^{\ell}
    + 4\alpha^2 p\right)
    +\sum_{k=1}^K4^{L_0}\alpha^2\\
    &\approx 4\alpha^2n_bKN+8\alpha^2p(L-L_0) + 4^{L_0}K\alpha^2.% = O(N).
  \end{aligned}
\end{equation}
As $L_0$ is small, thus the total number of parameters is $O(n_b\alpha^2KN)$.

%Here the number of parameters in $b$ in \eqref{eq:conv} is also ignored.

\paragraph{Translation-equivariant  case.}
For the linear system \eqref{eq:integral}, if the kernel is of convolution type, i.e.,
$a(x,y)=a(x-y)$, then $A$ is a cyclic matrix. So are the matrices $D_j\sps{\ell}$ and $A\sps{L}$.
Therefore, the $\LC$ layer in \cref{alg:wl} can be replaced with a $\Conv$ layer.

In the nonlinear case, the operator $A$ is {\em translation-equivariant} if
\begin{equation}\label{eq:invariant}
  \cT A (v) = A (\cT v)
\end{equation}
holds for any translation operator $\cT$. In this case, each $\LC$ layer in \cref{alg:bcrnet} is
replaced with a $\Conv$ layer with the same window size, while each $\Dense$ layer is replaced with
a $\Conv$ layer with window size equal to the input size. Since the number of parameters of a
$\Conv$ layer is $\alpha_{\out}\alpha_{\ipt}w$, the number of parameters of \cref{alg:bcrnet} is
\begin{equation}\label{eq:np_conv}
    \begin{aligned}
        N_{\mathrm{params}} &= \sum_{\ell=L_0}^{L-1}\left(4\alpha^2 p 
      + \sum_{k=1}^K(2\alpha)^2n_b
      + 4\alpha^2 p\right)
      +\sum_{k=1}^K4^{L_0}\alpha^2\\
      &\approx 4\alpha^2n_bK(L-L_0)+8\alpha^2p(L-L_0) + 4^{L_0}K\alpha^2\\
      & \approx O(n_b\alpha^2K\log(N)),% = O(\log N),
    \end{aligned}
\end{equation}
which is only logarithmic in $N$.

%==================================================================================================
\section{Applications}\label{sec:application}

%% : the diagonal
%% of the inverse matrix of an elliptic operator and the parametric $p$-Laplacian in homogenization
%% theory.
%All the tests are run on GPU with data type \verb|float32|.

We implement BCR-Net with Keras \cite{keras} (running on top of TensorFlow \cite{tensorflow}) using
Nadam as the optimizer \cite{dozat2015incorporating} and the mean squared error as the loss
function. The parameters in BCR-Net are initialized randomly from the normal distribution, and the
batch size is always set as two percent of the size of the training set.

In this section, we study the performance of BCR-Net using a few examples. In the experiments, the
support of the scaling function $\varphi(x)$ is chosen to be $[0,2p-1]$ with $p=3$. The activation
function in wavelet transform is set to be the identity as we mentioned, while ReLU is used in the
$\LC$ and $\Dense$ layers. As discussed in Section \ref{sec:bcrnet}, when the operator is
translation-equivariant, the $\LC$ layers and $\Dense$ layers in \cref{alg:bcrnet} are replaced by
$\Conv$ and fully connected $\Conv$ layers, respectively. The selection of parameters $\alpha$
(number of channels) and $K$ (number of $\LC$ layers in Algorithm \ref{alg:bcrnet}) are problem
dependent.

For each sample, let $u$ be the {\em exact} solution generated by numerical discretization of PDEs
and $u_{\NN}$ be the prediction from BCR-Net. Then the error of this specific sample is calculated
by the relative error measured in the $\ell^2$ norm:
\begin{equation}\label{eq:relativeerror}
  \epsilon = \frac{\|u-u_{\NN}\|_{\ell^2}}{\|u\|_{\ell^2}}.
\end{equation}
The training error $\trainerror$ and test error $\testerror$ are then obtained by averaging
\eqref{eq:relativeerror} over a given set of training or testing samples, respectively. The
numerical results presented in this section are obtained by repeating the training process five
times, using different random seeds.

%--------------
\subsection{Diagonal of the inverse matrix of elliptic operator}\label{sec:diag}

Elliptic operators of form $H=-\Delta + v(x)$ appear in many mathematical models. In quantum
mechanics, the Hamiltonian of the Schr\"odinger equation takes this form and it is fundamental in
describing the dynamics of quantum particles. In probability theory, this operator describes for
example the behavior of a random walk particle that interacts with the environment.

Here, we are interested the Green function $G=H^{-1}=(-\Delta + v(x))^{-1}$ and, more specifically,
the dependence of its diagonal $G(x,x)$ on the potential $v(x)$. The diagonal of the Green's
function plays an important role in several applications. For example, in density function theory,
when combined with appropriate rational expansions, the diagonal of the (shifted) Green's function
allows one to compute the Kohn-Sham map efficiently \cite{lin2009fast}. As another example, for
random walk particles that decay at a spatially dependent rate given by $v(x)>0$, the diagonal
$G(x,x)$ gives the expected local time (or the number of visits in the discrete setting) at position
$x$ for a particle starting from $x$, via the Feynman-Kac formula \cite{evans2013}.

In the following numerical studies, we work with the operator $H = -\Delta + v(x)$ on $[0,1]^2$ in
2D with the periodic boundary condition. The differential operator is discretized using the standard
five-point central difference with $80$ points per dimension. The potential $v$ is generated by
randomly sampling on a $10\times 10$ grid independently from the standard Gaussian distribution
$\cN(0,1)$, interpolating it to the $80\times 80$ grid via Fourier interpolation, and finally
performing a pointwise exponentiation.

By denoting the diagonal of $G=H^{-1}$ by $g(x)$, we apply BCR-Net to learn the nonlinear map from
the potential $v(x)$ to $g(x)$:
\begin{equation}\label{eq:diagonal_map}
  v(x) \rightarrow g(x):=G(x,x).
\end{equation}
The approximation to $g(x)$ produced by BCR-Net will be denoted by $g_{\NN}(x)$.

%\cref{fig:greenSolX} presents a sample of the potential.

\begin{figure}[htb]
    \centering
    \subfloat[\label{fig:greenAl} $K=5$]{
    \begin{overpic}[width=0.4\textwidth, clip]{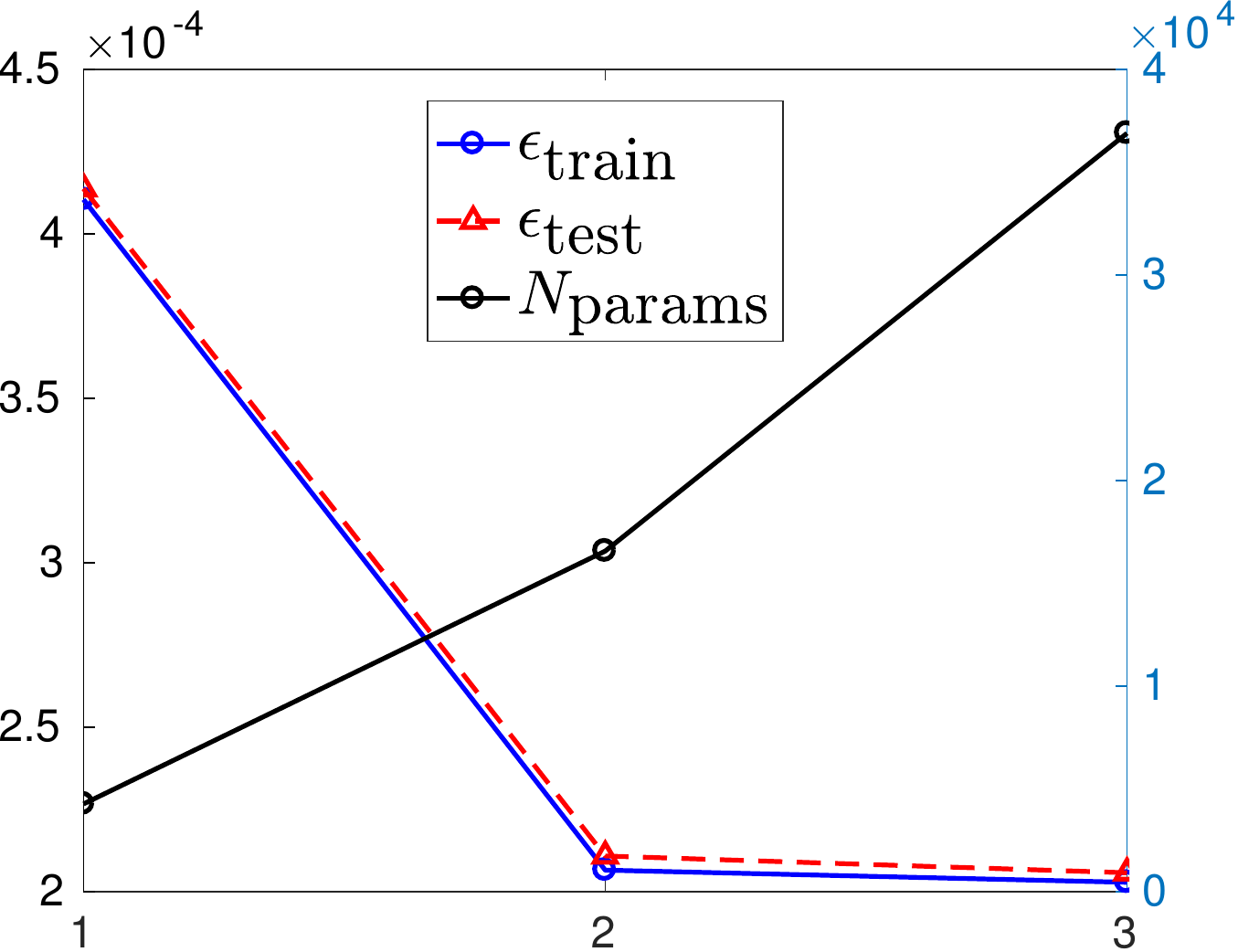}
        \put(60,-1){$\alpha$}
        \put(-2,53){$\epsilon$}
        \put(94,48){$N_{\mathrm{params}}$}
    \end{overpic}
    }
    \qquad\qquad
    \subfloat[\label{fig:greenK} $\alpha=2$]{
    \begin{overpic}[width=0.4\textwidth, clip]{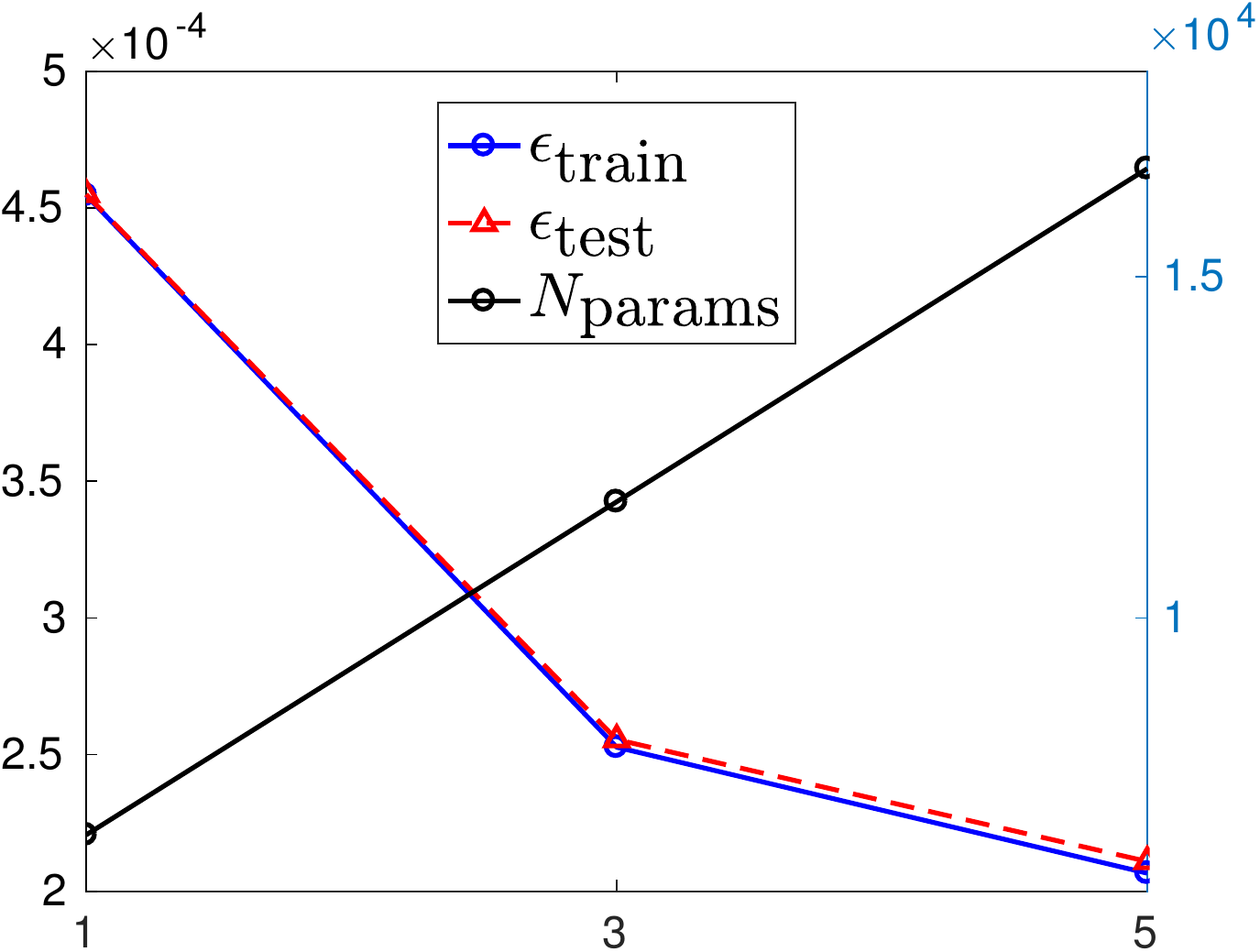}
        \put(60,-1){$K$}
        \put(-2,53){$\epsilon$}
        \put(94,48){$N_{\mathrm{params}}$}
    \end{overpic}
    }
    \caption{\label{fig:greenAlK} The relative error and the number of parameters of BCR-Net for
      \eqref{eq:diagonal_map} with $\Ntrainsample = \Ntestsample =20000$. (a) Different channel
      numbers ($\alpha$). (b) Different $\LC$ layer numbers ($K$).}
\end{figure}

%%%%%%% choose one from the following two figures %%%%%%%%%%%%%%
%%%%% 3D figure, plotted by ``mesh'' in matlab
\begin{figure}[ht]
    \centering
    \subfloat[\label{fig:greenSolX}$v(x)$]{
      \begin{overpic}[clip, width=0.32\textwidth]{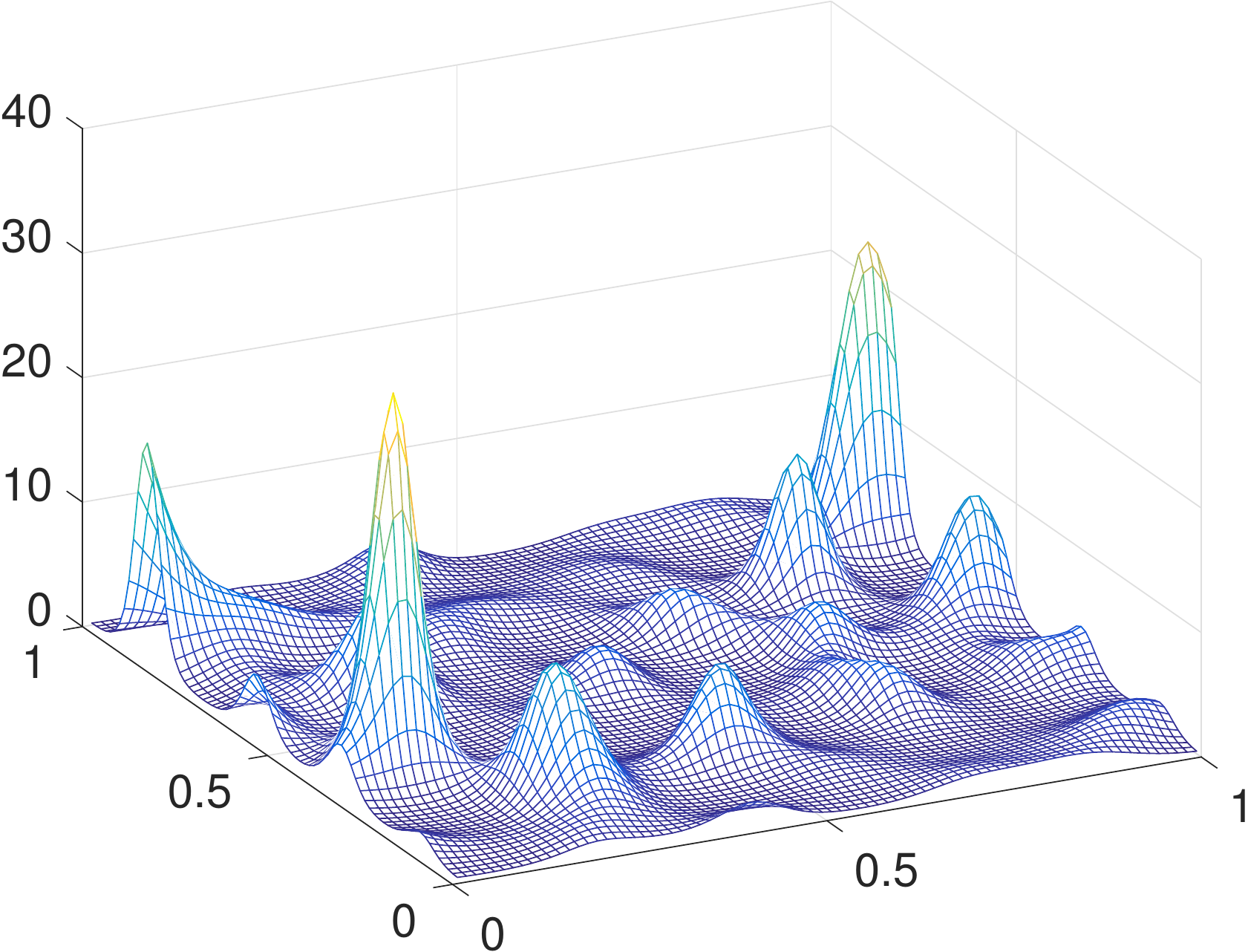}
      \end{overpic}
    }
    \subfloat[$g_{\NN}(x)$]{
      \begin{overpic}[clip, width=0.32\textwidth]{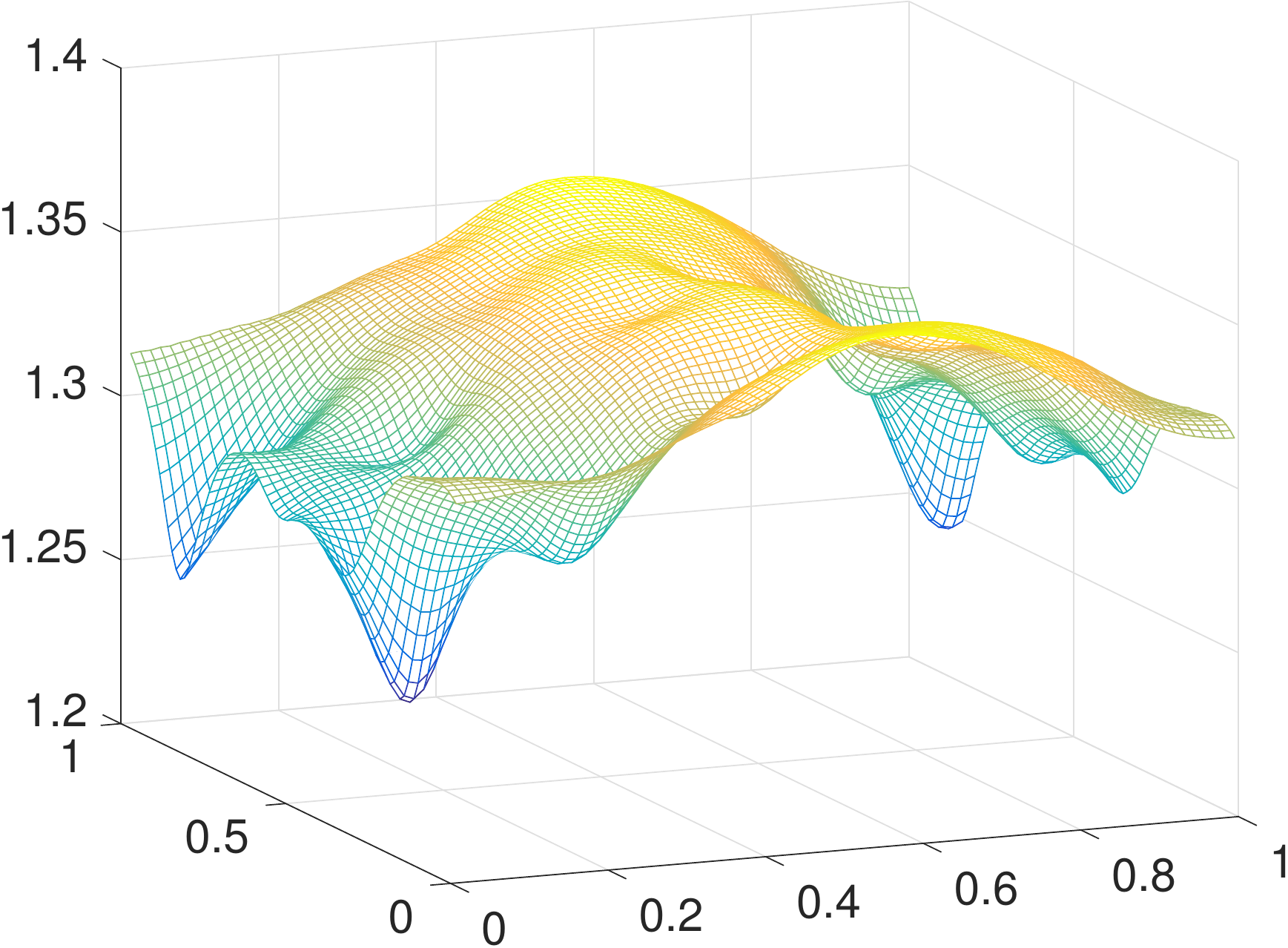}
      \end{overpic}
    }
    \subfloat[$g(x)-g_{\NN}(x)$]{
      \begin{overpic}[clip, width=0.32\textwidth]{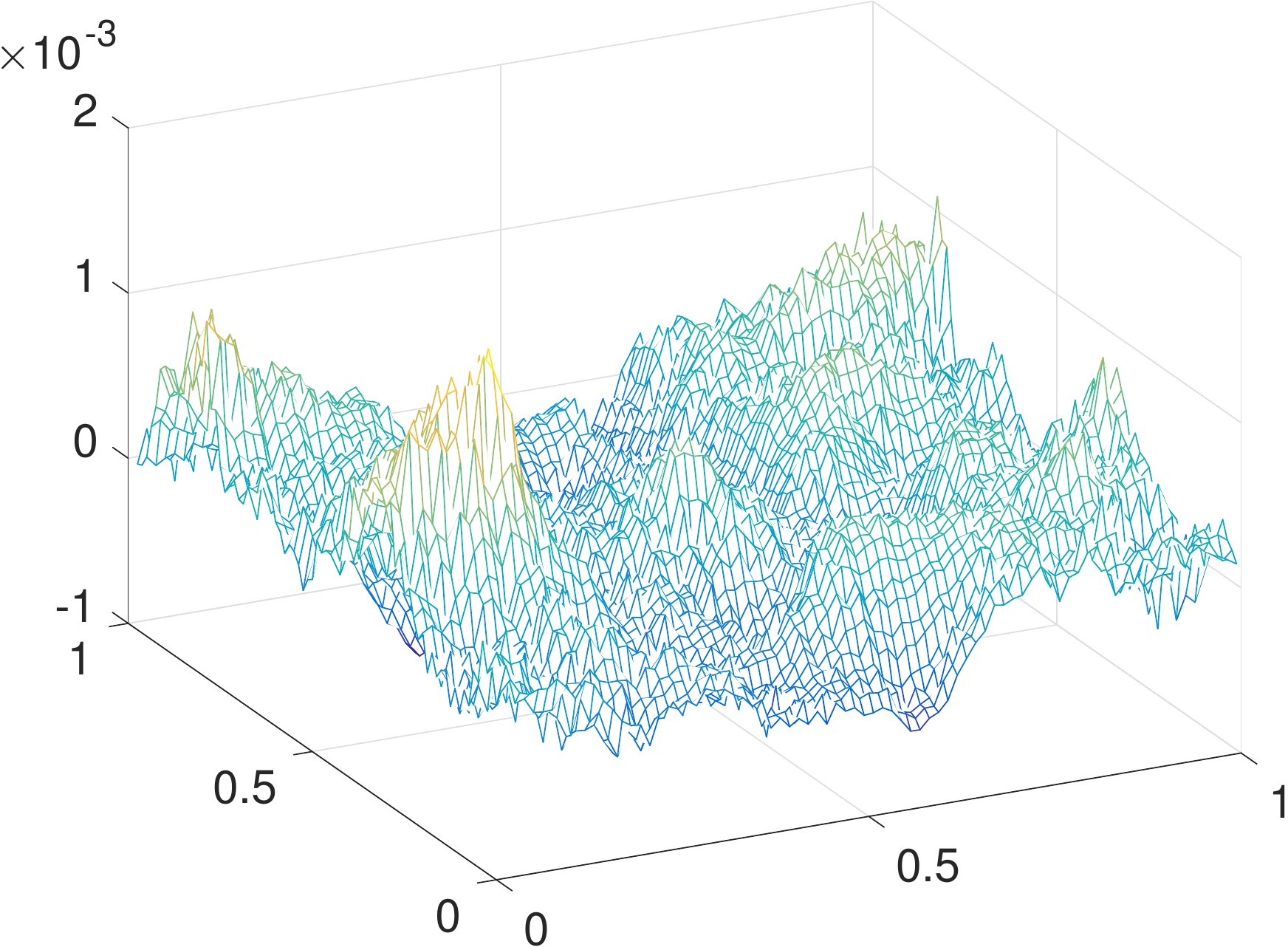}
      \end{overpic}
    }
    \caption{\label{fig:greenSol} A random sample potential $v(x)$, the prediction $g_{\NN}(x)$ from
      BCR-Net with $\alpha=2$ and $K=5$, and the error from the reference solution of
      \eqref{eq:diagonal_map}.}
\end{figure}

\cref{fig:greenAlK} presents the numerical results with different choices for $\alpha$ (channel
number) and $K$ (layer number) in \cref{alg:bcrnet}. The band of the $D\sps{\ell}_j$ matrices is set
to be $n_b=3$. When $\alpha$ or $K$ increases, the error decreases consistently. \cref{fig:greenAl}
plots the results for different $\alpha$ values with $K$ fixed. Note that $\alpha=2$ already
achieves a fairly accurate result.  \cref{fig:greenK} shows the results for different $K$ values with
$\alpha$ fixed, showing that deeper networks clearly give better results. In each case, the test
error is close to the training error, implying that there is no over-fitting in our model. As
discussed in \cref{sec:bcrnet}, the number of parameters is proportional to $\alpha^2K$, which
agrees with the curves in \cref{fig:greenAlK}.  In \cref{fig:greenSol}, we plot a sample of the
potential $v(x)$, along with the prediction from BCR-Net and its error in comparison with the
reference solution ($\alpha=2$ and $K=5$).

It is worth pointing out that the test error is quite small (around $2.1\times 10^{-4}$ with
$\alpha=2$ and $K=5$), while the number of parameters is $\Nparams=1.6\times 10^{4}$.  Comparing
with millions of parameters in the applications on images \cite{Krizhevsky2012,
  Ronneberger2015}, BCR-Net only uses tens of thousands of parameters.

%% \fy{to lexing: it is really millions. 
%% The network in \cite{Krizhevsky2012} is globally connected. 
%% I have calculated the number of parameters of U-net in \cite{Ronneberger2015}, \#params is $2.3\times
%% 10^7$. For a lot of famous network, \#params are all $10^6\sim 10^8$.}

%--------------
\subsection{Nonlinear homogenization theory} 

The second example is concerned with homogenization theory, which studies effective models for
differential and integral equations with oscillatory coefficients. In the simplest setting, consider
the linear second order elliptic PDE in a domain $\Omega$ 
\[
-\nabla\cdot \left(a\left(\frac{x}{\eps}\right) \nabla u^\eps(x) \right) = 0
\]
with appropriate boundary conditions on $\partial\Omega$, where $a(\cdot)$ is a periodic function on
the unit cube $[0,1]^d$. The homogenization theory
\cite{Bensoussan-2011,Engquist-2008,Jikov-1994,Pavliotis-2008} states that, when $\eps$ goes to zero,
the solution $u^\eps(x)$ exhibits a multiscale decomposition
\[
u^\eps(x) = u_0(x) + \eps u_1(x) + \eps^2 u_2(x) + \cdots.
\]
Here $u_0(x)$ is the solution of a constant elliptic PDE
\[
-\nabla\cdot (A_0 \nabla u_0(x)) = 0
\]
where the {\em constant} matrix $A_0$ is called the {\em effective coefficient tensor}. The next
term $u_1(x)$ that depends on the gradient of $u_0(x)$ can be written as
\[
u_1(x)= \sum_{i=1}^d \eta_i\left(\frac{x}{\eps}\right) \nabla_i u_0(x)
\]
in terms of the so-called {\em corrector functions} $\{\eta_i(x)\}_{i=1,\ldots,d}$.  For each
$i\in\{1,\ldots,d\}$, the corrector function $\eta_i(x)$ is the solution of the periodic problem
\begin{equation}\label{eq:corrector}
  -\nabla\cdot (a(x) (\nabla \eta_i(x) + e_i)) = 0, \quad \int_{[0,1]^d}\eta_i(x)\dd x = 0,
\end{equation}
where $e_i$ is the canonical basis vector in the $i$-th coordinate. Both the effective coefficient
tensor $A_0$ and the correctors $\eta_i(x)$ are important for studying multiscale problems in
engineering applications. Note from \eqref{eq:corrector} that computing the corrector $\eta_i(x)$
for each $i$ requires a single PDE solve over the periodic cube.

For many nonlinear problems, a similar homogenization theory holds. Consider for
example the variational problem
\[
u_\eps = \argmin_v \int_\Omega f\left( \frac{x}{\eps}, \nabla v \right) \dd x
\]
with appropriate boundary conditions on $\partial\Omega$, where $f\left( \frac{x}{\eps}, \nabla v
\right)$ is given by
\[
f\left( \frac{x}{\eps}, \nabla u \right) = a \left(\frac{x}{\eps}\right) |\nabla u|^p
\]
for some constant $p > 1$. In this case, the corrector function is parameterized by a unit vector
$\xi \in \mathbb{S}^{d-1}$: for a fixed unit vector $\xi$, the corrector function $\chi_{p,\xi}(x)$
satisfies
%% is concerned with the derivation of equations for average of solutions of equations with rapidly
%% varying coefficients, which arises in obtaining macroscopic equations for systems with a fine
%% microscopic structure \cite{sanchez1980non}.  Here we study the homogenization of $p$-Laplacian with
%% the form
\begin{equation}\label{eq:elliptic}
  -\nabla\cdot \left( a(x)|\nabla \chi_{p,\xi}(x)+\xi|^{p-2}(\nabla \chi_{p,\xi}(x)+\xi) \right) = 0,\quad
  \int_{[0,1)^d}\chi_{p,\xi}(x)\dd x = 0,
\end{equation}
with the periodic boundary condition. Since \eqref{eq:elliptic} is a nonlinear PDE, the numerical
solution of the corrector function in the nonlinear case is computationally more challenging.
%% In practical applications, $a(x)$ can be uncertainty quantity. Therefore \eqref{eq:elliptic} is
%% solved numerically an exponential number of times.

%In the numerical tests, $p=3$ and $\xi$ is chosen to be the unit vector in the $x_1$ direction.

Here, we apply BCR-Net to learn the map from the periodic oscillatory coefficient function $a(x)$ to
the corrector $\chi_{p,\xi}(x)$ for a given $\xi=(\xi_1,\ldots,\xi_d)^T$. Note that, for the case
$p=2$, \eqref{eq:elliptic} reduced to \eqref{eq:corrector}. It is noticed that, though quite
different numerically, the solution of the linear problem \eqref{eq:corrector} serves as a
reasonable baseline for the nonlinear problem \eqref{eq:elliptic}. Motivated by this observation, we
compute $\eta_{\xi}(x):=\sum_{i=1}^d\xi_i\eta_i(x)$ by solving the linear system efficiently and use
BCR-Net only to learn the map from $a(x)$ to the difference $g(x)$ between the nonlinear and linear
correctors
\begin{equation}\label{eq:elliptic_map}
  a(x)\rightarrow g(x) := \chi_{p,\xi}(x) - \eta_{\xi}(x) = \chi_{p,\xi}(x)-\sum_{i=1}^d\xi_i\eta_i(x).
\end{equation}

\paragraph{Two-dimensional case.}
In this test, we set $p=3$ and discretize the coefficient field $a(x)$ with an $80\times 80$
Cartesian grid. Each realization of $a(x)$ is generated by randomly sampling on a $10\times 10$ grid
with respect to $\cN(0,1)$, interpolating it to a $80\times 80$ grid via the Fourier transform, and
finally taking pointwise exponential.

\begin{figure}[ht]
    \centering
    \subfloat[\label{fig:homoAl} $K=5$]{
     \begin{overpic}[width=0.4\textwidth, clip]{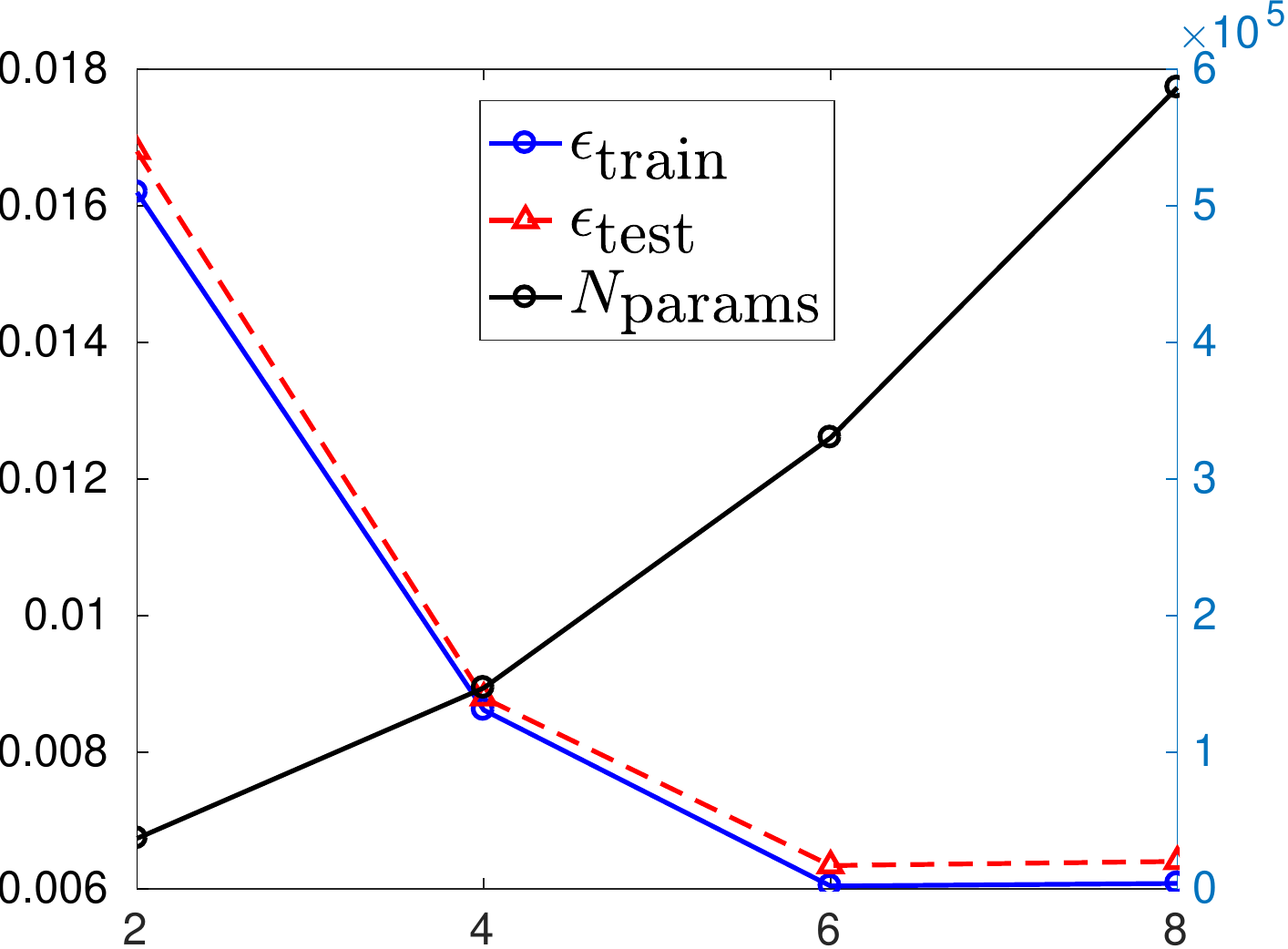}
         \put(50,-1){$\alpha$}
         \put(-2,53){$\epsilon$}
         \put(98,48){$N_{\mathrm{params}}$}
     \end{overpic}
    }
    \qquad\qquad
    \subfloat[\label{fig:homoK} $\alpha=6$]{
    \begin{overpic}[width=0.4\textwidth, clip]{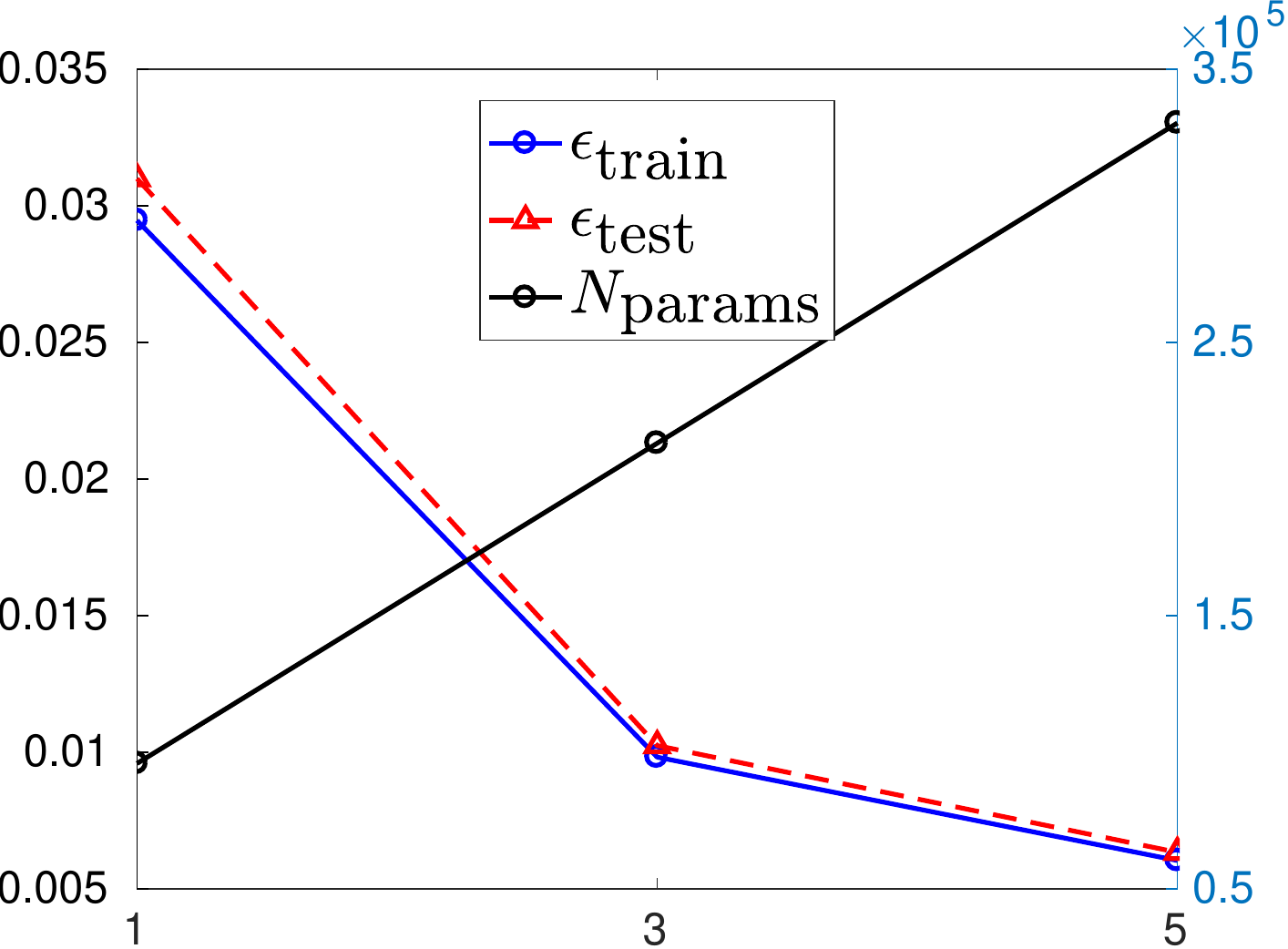}
        \put(60,-1){$K$}
        \put(-2,53){$\epsilon$}
        \put(98,48){$N_{\mathrm{params}}$}
    \end{overpic}
    }
    \caption{\label{fig:homoAlK} Relative error and number of parameters of BCR-Net for
      \eqref{eq:elliptic_map} with $\Ntrainsample = \Ntestsample =20000$ in 2D.  (a) Different
      channel numbers ($\alpha$). (b) Different $\LC$ layers numbers ($K$).}
\end{figure}

\begin{figure}[ht]
    \centering
    \subfloat[\label{fig:homoSolX}$a(x)$]{
      \begin{overpic}[clip, width=0.32\textwidth]{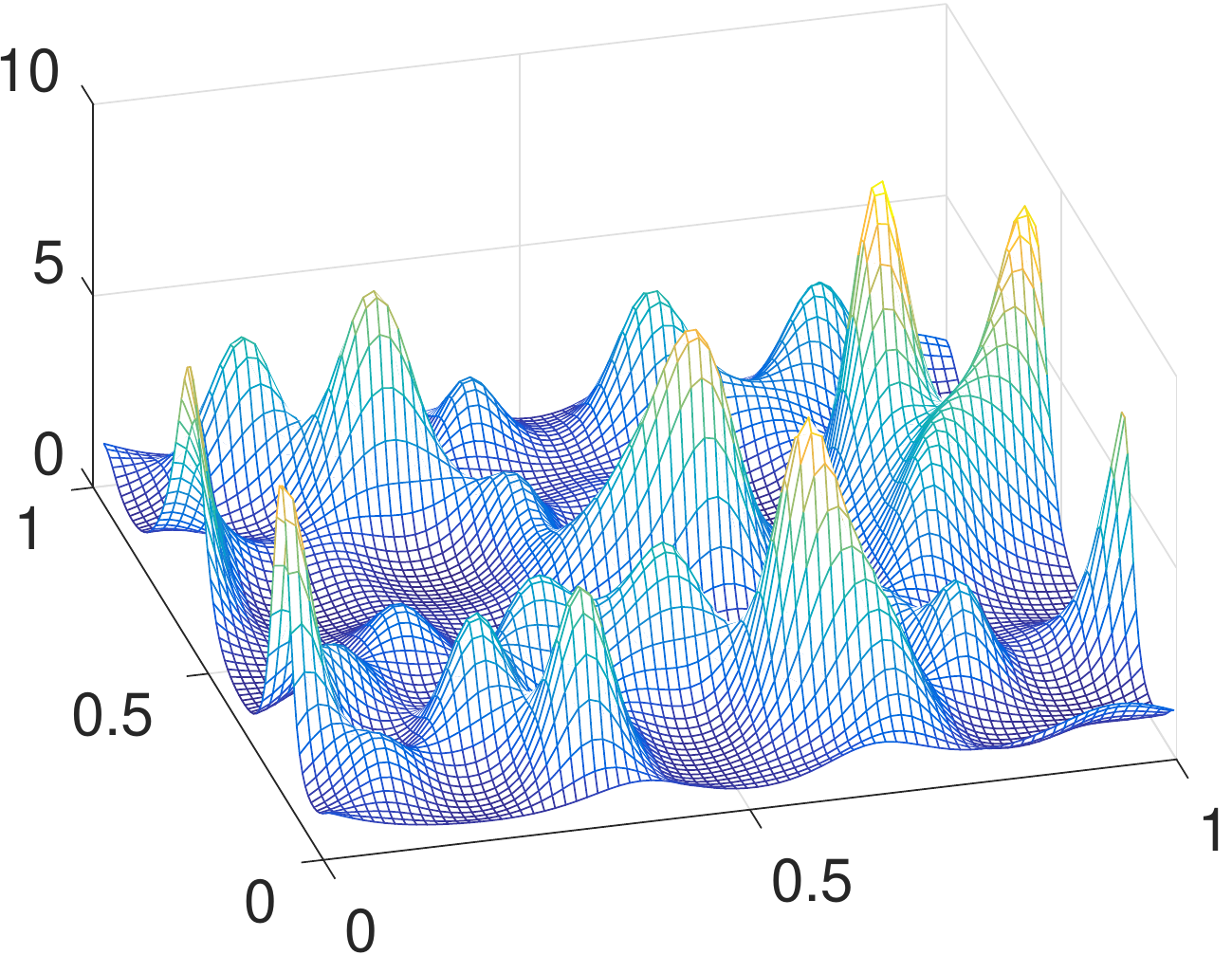}
      \end{overpic}
    }
    \subfloat[$\chi_{\xi,p,\NN}(x)$]{
      \begin{overpic}[clip, width=0.32\textwidth]{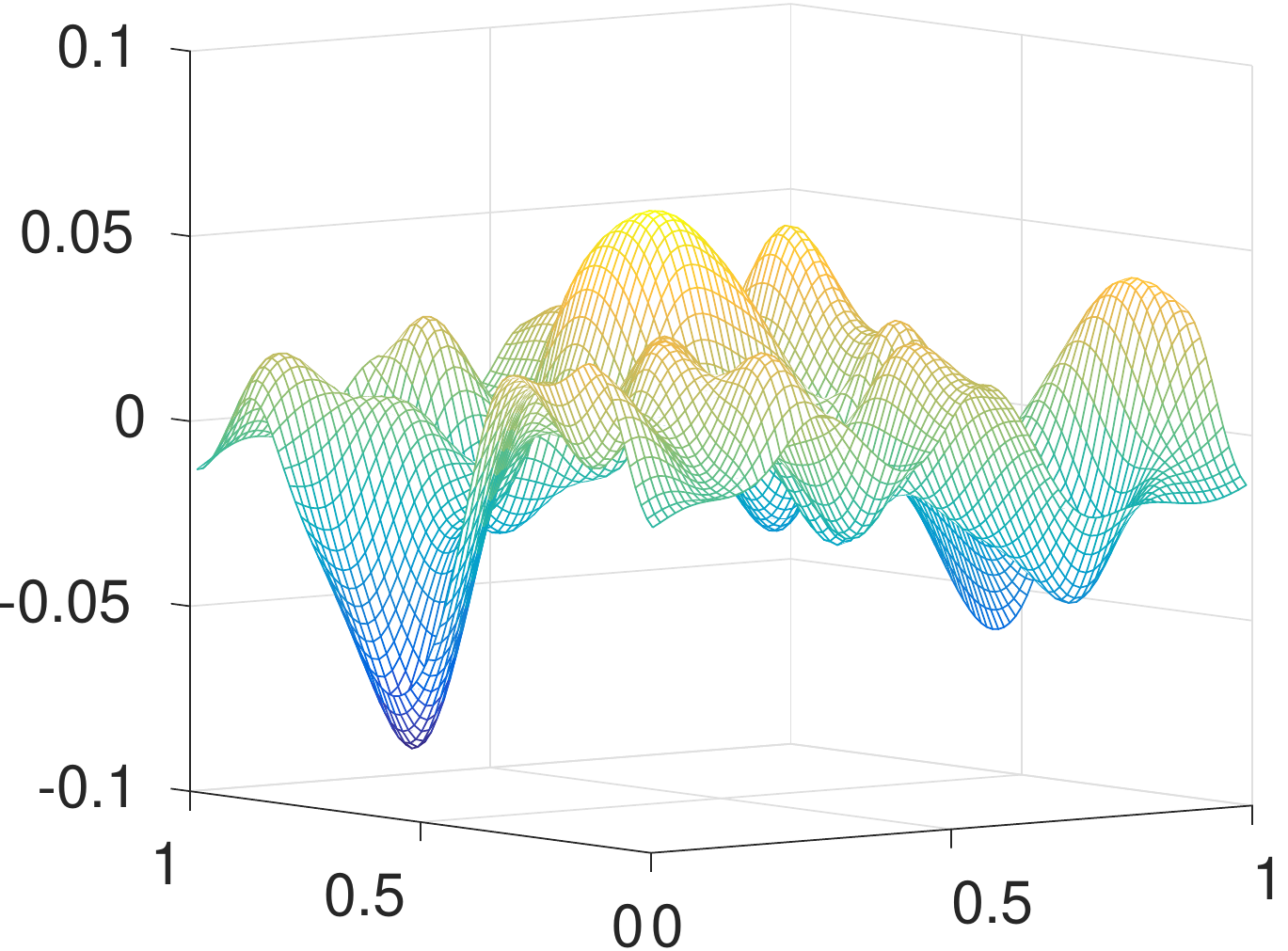}
      \end{overpic}
    }
    \subfloat[$\chi_{\xi,p}(x)-\chi_{\xi,p,\NN}(x)$]{
      \begin{overpic}[clip, width=0.32\textwidth]{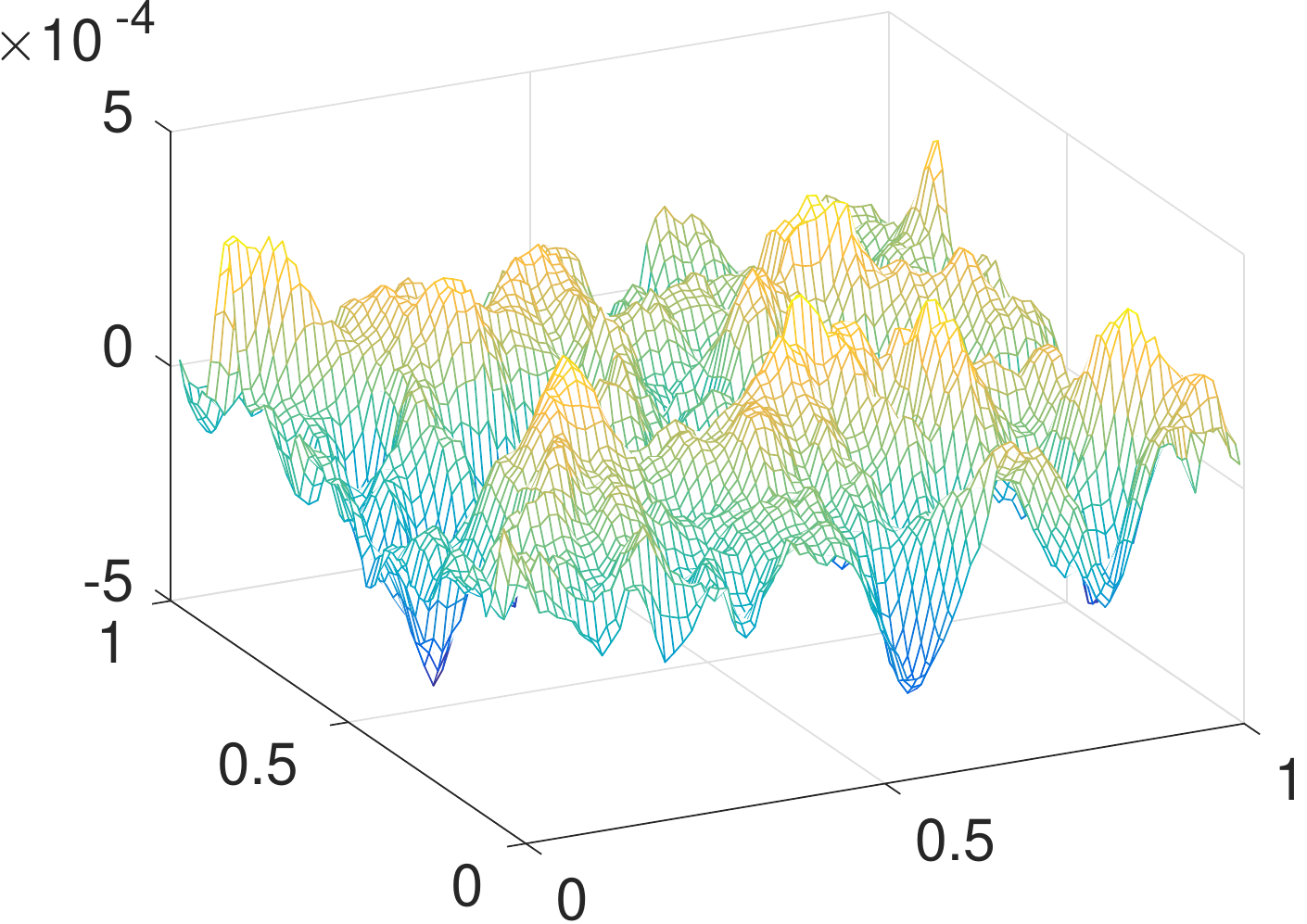}
      \end{overpic}
    }
    \caption{\label{fig:homoSol} A random coefficient $a(x)$, the approximation
      $\chi_{\xi,p,\NN}(x)$ produced by BCR-Net with $\alpha=6$ and $K=5$, and the error from the
      the reference solution of \eqref{eq:elliptic_map} in 2D.}
\end{figure}

\cref{fig:homoAlK} summarizes the numerical results for different choices of $\alpha$ and $K$ in
\cref{alg:bcrnet} with $n_b=5$. Notice that the error behavior is comparable to the one shown in
\cref{sec:diag}. The relative error decreases consistently as $\alpha$ or $K$ increases. In
addition, there is no sign of over-fitting and the number of parameters grows proportional to
$\alpha^2K$ in agreement with the complexity analysis. From \cref{fig:homoAlK}, we notice that the
best choice of the parameters for this problem is $\alpha=6$ and $K=5$. \cref{fig:homoSol} shows a
random sample of the coefficient field $a(x)$, the approximation to $\chi_{p,\xi}(x)$ predicted by
BCR-Net, and the error when compared with the reference solution.

\begin{figure}[htbp]
    \centering \subfloat[\label{fig:homo3DSolX}$a(x)$]{
      \begin{overpic}[clip, width=0.32\textwidth]{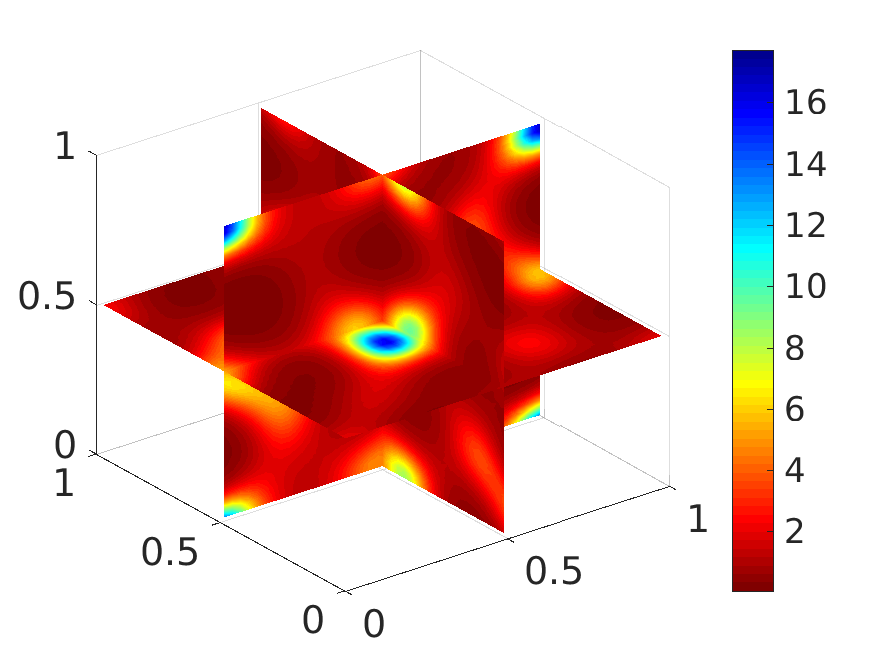}
    \end{overpic}
    }
    \subfloat[$\chi_{\xi,p,\NN}(x)$]{
      \begin{overpic}[clip, width=0.32\textwidth]{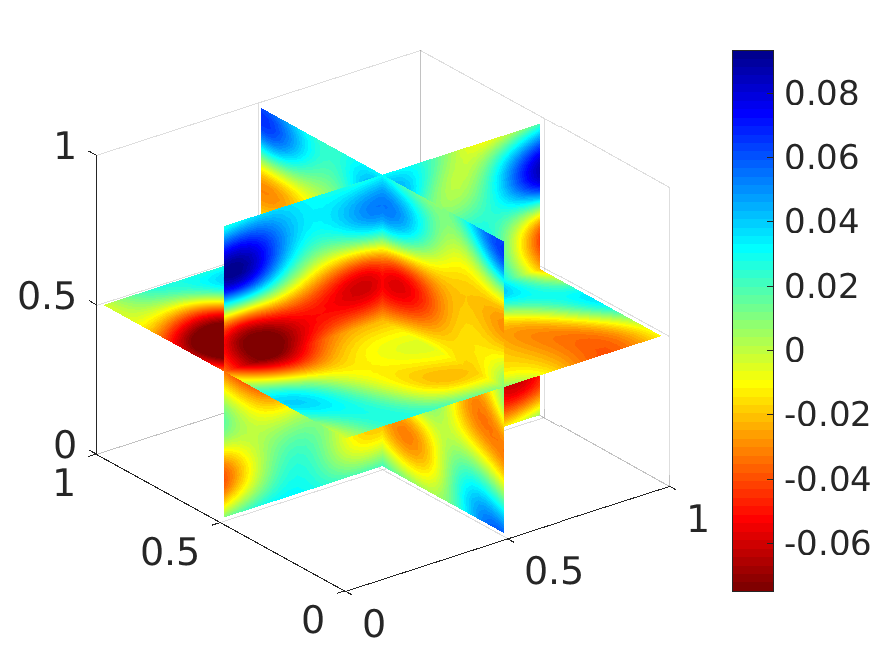}
      \end{overpic}
    }
    \subfloat[$\chi_{\xi,p}(x)$]{
      \begin{overpic}[clip, width=0.32\textwidth]{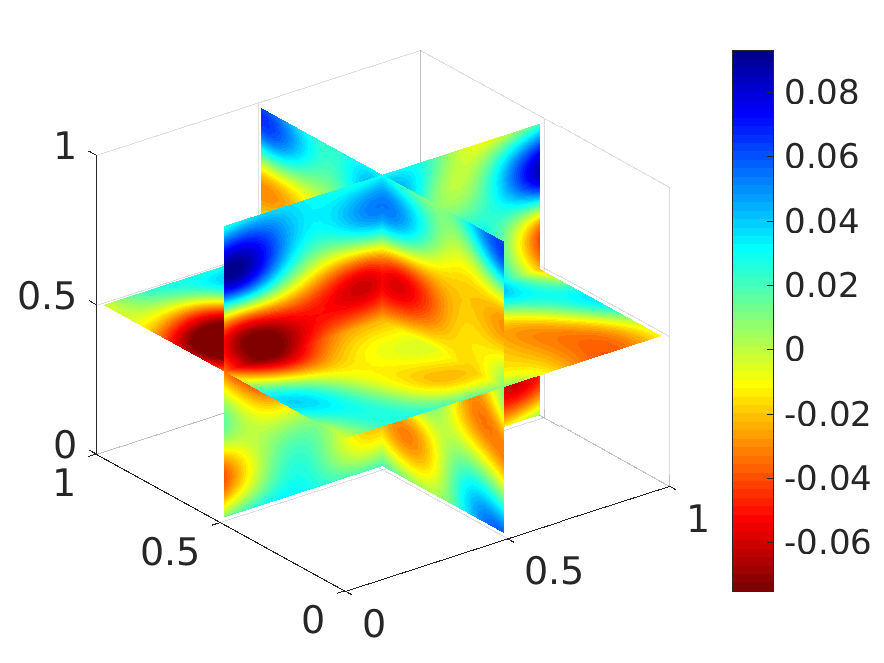}
      \end{overpic}
    }
    \caption{\label{fig:homo3DSol}
      A random coefficient $a(x)$ in 3D, the prediction $\chi_{\xi,p,\NN}(x)$ from BCR-Net with
      $\alpha=4$ and $K=5$, and the reference solution of \eqref{eq:elliptic_map}.}
\end{figure}

\paragraph{Three-dimensional case.}
We set the discretization grid for $a(x)$ to be $40\times 40\times 40$ and $p$ to be 3. The
coefficient field $a(x)$ is generated by randomly sampling on a $5\times5\times5$ grid from
$\cN(0,1)$, interpolating to the full grid via Fourier transformation, and finally taking pointwise
exponential. The BCR-Net is trained with $n_b=3$, $\alpha=4$, and $K=5$ using
$\Ntrainsample=\Ntestsample=10000$ samples. The number of parameters is $5.7\times 10^{5}$ and the
test error is $8.7\times 10^{-3}$.  As \cref{fig:homo3DSol} shows the results for one random
realization of $a(x)$, BCR-Net is able to reproduce the corrector function quite accurately.

%==================================================================================================
\section{Conclusion}\label{sec:conclusion}

In this paper, inspired by the nonstandard wavelet form proposed by Beylkin, Coifman, and Rokhlin in
\cite{bcr}, we developed a novel neural network architecture, called BCR-Net, to approximate certain
nonlinear generalization of pseudo-differential operators. This NN demonstrates promising results
while approximating the nonlinear maps arising from homogenization theory and stochastic
computation, using only tens of thousands of parameters.

The BCR-Net architecture can be naturally extended in several ways. For instance, the $\LC$ layers
can be altered to include other network structures, such as parallel sub-networks or the ResNet
architecture \cite{he2016deep}. The $\Conv$ layers corresponding the wavelet transforms can also be
replaced with other types of building blocks.

%==================================================================================================
\section*{Acknowledgments}
The work of Y.F. and L.Y. is partially supported by the U.S. Department of Energy, Office of
Science, Office of Advanced Scientific Computing Research, Scientific Discovery through Advanced
Computing (SciDAC) program.  The work of L.Y. is also partially supported by the National Science
Foundation under award DMS-1818449. The work of C. O. is partially supported by the Stanford
Graduate Fellowship in Science \& Engineering.  This work is also supported by the GCP Research
Credits Program from Google.

\bibliographystyle{abbrv}
\bibliography{nn}

\def\cprime{$'$}
\begin{thebibliography}{10}

\bibitem{tensorflow}
M.~Abadi et~al.
\newblock Tensorflow: A system for large-scale machine learning.
\newblock In {\em OSDI}, volume~16, pages 265--283, 2016.

\bibitem{alpert2002adaptive}
B.~Alpert, G.~Beylkin, D.~Gines, and L.~Vozovoi.
\newblock Adaptive solution of partial differential equations in multiwavelet
  bases.
\newblock {\em Journal of Computational Physics}, 182(1):149--190, 2002.

\bibitem{Araya-Polo2018}
M.~Araya-Polo, J.~Jennings, A.~Adler, and T.~Dahlke.
\newblock Deep-learning tomography.
\newblock {\em The Leading Edge}, 37(1):58--66, 2018.

\bibitem{badrinarayanan2015segnet}
V.~Badrinarayanan, A.~Kendall, and R.~Cipolla.
\newblock Segnet: A deep convolutional encoder-decoder architecture for image
  segmentation.
\newblock {\em IEEE Transactions on Pattern Analysis and Machine Intelligence},
  2017.

\bibitem{Bensoussan-2011}
A.~Bensoussan, J.-L. Lions, and G.~Papanicolaou.
\newblock {\em Asymptotic analysis for periodic structures}.
\newblock AMS Chelsea Publishing, Providence, RI, 2011.
\newblock Corrected reprint of the 1978 original [MR0503330].

\bibitem{berg2017unified}
J.~Berg and K.~Nystr{\"o}m.
\newblock A unified deep artificial neural network approach to partial
  differential equations in complex geometries.
\newblock {\em Neurocomputing}, 317:28--41, 2018.

\bibitem{bcr}
G.~Beylkin, R.~Coifman, and V.~Rokhlin.
\newblock Fast wavelet transforms and numerical algorithms {I}.
\newblock {\em Communications on pure and applied mathematics}, 44(2):141--183,
  1991.

\bibitem{Bruna2012}
J.~Bruna and S.~Mallat.
\newblock Invariant scattering convolution networks.
\newblock {\em IEEE Transactions on Pattern Analysis and Machine Intelligence},
  35(8):1872--1886, 2013.

\bibitem{Chen2018DeepLab}
L.~C. Chen, G.~Papandreou, I.~Kokkinos, K.~Murphy, and A.~L. Yuille.
\newblock Deeplab: Semantic image segmentation with deep convolutional nets,
  atrous convolution, and fully connected crfs.
\newblock {\em IEEE Transactions on Pattern Analysis and Machine Intelligence},
  40(4):834--848, 2018.

\bibitem{keras}
F.~Chollet et~al.
\newblock Keras.
\newblock \url{https://keras.io}, 2015.

\bibitem{cohen2003numerical}
A.~Cohen.
\newblock {\em Numerical analysis of wavelet methods}, volume~32.
\newblock Elsevier, 2003.

\bibitem{CohenSharir2018}
N.~Cohen, O.~Sharir, and A.~Shashua.
\newblock On the expressive power of deep learning: A tensor analysis.
\newblock {\em arXiv preprint arXiv:1603.00988}, 2018.

\bibitem{daubechies1988orthonormal}
I.~Daubechies.
\newblock Orthonormal bases of compactly supported wavelets.
\newblock {\em Communications on pure and applied mathematics}, 41(7):909--996,
  1988.

\bibitem{dozat2015incorporating}
T.~Dozat.
\newblock Incorporating nesterov momentum into adam.
\newblock {\em International Conference on Learning Representations}, 2016.

\bibitem{Engquist-2008}
B.~Engquist and P.~E. Souganidis.
\newblock Asymptotic and numerical homogenization.
\newblock {\em Acta Numer.}, 17:147--190, 2008.

\bibitem{evans2013}
L.~C. Evans.
\newblock {\em An introduction to stochastic differential equations}.
\newblock American Mathematical Society, Providence, RI, 2013.

\bibitem{fan2018mnnh2}
Y.~Fan, J.~Feliu-Fab{\`a}, L.~Lin, L.~Ying, and L.~Zepeda-N{\'u}\~nez.
\newblock A multiscale neural network based on hierarchical nested bases.
\newblock {\em arXiv preprint arXiv:1808.02376}, 2018.

\bibitem{fan2018mnn}
Y.~Fan, L.~Lin, L.~Ying, and L.~Zepeda-N\'u\~nez.
\newblock A multiscale neural network based on hierarchical matrices.
\newblock {\em arXiv preprint arXiv:1807.01883}, 2018.

\bibitem{goodfellow2016deep}
I.~Goodfellow, Y.~Bengio, A.~Courville, and Y.~Bengio.
\newblock {\em Deep learning}, volume~1.
\newblock MIT press Cambridge, 2016.

\bibitem{han2018solving}
J.~Han, A.~Jentzen, and W.~E.
\newblock Solving high-dimensional partial differential equations using deep
  learning.
\newblock {\em Proceedings of the National Academy of Sciences},
  115(34):8505--8510, 2018.

\bibitem{he2016deep}
K.~He, X.~Zhang, S.~Ren, and J.~Sun.
\newblock Deep residual learning for image recognition.
\newblock In {\em Proceedings of the IEEE conference on computer vision and
  pattern recognition}, pages 770--778, 2016.

\bibitem{Hinton2012}
G.~Hinton, L.~Deng, D.~Yu, G.~E. Dahl, A.~r.~Mohamed, N.~Jaitly, A.~Senior,
  V.~Vanhoucke, P.~Nguyen, T.~N. Sainath, and B.~Kingsbury.
\newblock Deep neural networks for acoustic modeling in speech recognition: The
  shared views of four research groups.
\newblock {\em IEEE Signal Processing Magazine}, 29(6):82--97, 2012.

\bibitem{Hornik91}
K.~Hornik.
\newblock Approximation capabilities of multilayer feedforward networks.
\newblock {\em Neural Networks}, 4(2):251--257, 1991.

\bibitem{Jikov-1994}
V.~V. Jikov, S.~M. Kozlov, and O.~A. Oleinik.
\newblock {\em Homogenization of differential operators and integral
  functionals}.
\newblock Springer-Verlag, Berlin, 1994.
\newblock Translated from the Russian by G. A. Yosifian [G. A.
  Iosif{\cprime}yan].

\bibitem{khoo2017solving}
Y.~Khoo, J.~Lu, and L.~Ying.
\newblock Solving parametric {PDE} problems with artificial neural networks.
\newblock {\em arXiv preprint arXiv:1707.03351}, 2017.

\bibitem{Khrulkov2018}
V.~Khrulkov, A.~Novikov, and I.~Oseledets.
\newblock Expressive power of recurrent neural networks.
\newblock {\em arXiv:1711.00811}, 2018.

\bibitem{Krizhevsky2012}
A.~Krizhevsky, I.~Sutskever, and G.~E. Hinton.
\newblock Imagenet classification with deep convolutional neural networks.
\newblock In {\em Proceedings of the 25th International Conference on Neural
  Information Processing Systems - Volume 1}, NIPS'12, pages 1097--1105, USA,
  2012. Curran Associates Inc.

\bibitem{leCunn2015}
Y.~LeCun, Y.~Bengio, and G.~Hinton.
\newblock Deep learning.
\newblock {\em Nature}, 521(436), 2015.

\bibitem{Leung2014}
M.~K.~K. Leung, H.~Y. Xiong, L.~J. Lee, and B.~J. Frey.
\newblock Deep learning of the tissue-regulated splicing code.
\newblock {\em Bioinformatics}, 30(12):i121--i129, 2014.

\bibitem{Yingzhou2018}
Y.~Li, X.~Cheng, and J.~Lu.
\newblock Butterfly-{Net}: Optimal function representation based on
  convolutional neural networks.
\newblock {\em arXiv preprint arXiv:1805.07451}, 2018.

\bibitem{lin2009fast}
L.~Lin, J.~Lu, L.~Ying, R.~Car, and W.~E.
\newblock Fast algorithm for extracting the diagonal of the inverse matrix with
  application to the electronic structure analysis of metallic systems.
\newblock {\em Communications in Mathematical Sciences}, 7(3):755--777, 2009.

\bibitem{LITJENS201760}
G.~Litjens, T.~Kooi, B.~E. Bejnordi, A.~A.~A. Setio, F.~Ciompi, M.~Ghafoorian,
  J.~A.~W.~M. van~der Laak, B.~van Ginneken, and C.~I. S\'anchez.
\newblock A survey on deep learning in medical image analysis.
\newblock {\em Medical Image Analysis}, 42:60--88, 2017.

\bibitem{MaSheridan2015}
J.~Ma, R.~P. Sheridan, A.~Liaw, G.~E. Dahl, and V.~Svetnik.
\newblock Deep neural nets as a method for quantitative structure–activity
  relationships.
\newblock {\em Journal of Chemical Information and Modeling}, 55(2):263--274,
  2015.

\bibitem{Mallat2009}
S.~Mallat.
\newblock {\em A wavelet tour of signal processing: the sparse way}.
\newblock Academic press, Boston, third edition, 2008.

\bibitem{Mhaskar2018}
H.~Mhaskar, Q.~Liao, and T.~Poggio.
\newblock Learning functions: When is deep better than shallow.
\newblock {\em arXiv preprint arXiv:1603.00988}, 2018.

\bibitem{misiti2013wavelets}
M.~Misiti, Y.~Misiti, G.~Oppenheim, and J.-M. Poggi.
\newblock {\em Wavelets and their Applications}.
\newblock John Wiley \& Sons, 2013.

\bibitem{Pavliotis-2008}
G.~A. Pavliotis and A.~M. Stuart.
\newblock {\em Multiscale methods}, volume~53 of {\em Texts in Applied
  Mathematics}.
\newblock Springer, New York, 2008.
\newblock Averaging and homogenization.

\bibitem{Raissi2018}
M.~Raissi and G.~E. Karniadakis.
\newblock Hidden physics models: Machine learning of nonlinear partial
  differential equations.
\newblock {\em Journal of Computational Physics}, 357:125 -- 141, 2018.

\bibitem{Ronneberger2015}
O.~Ronneberger, P.~Fischer, and T.~Brox.
\newblock U-net: Convolutional networks for biomedical image segmentation.
\newblock In {\em Medical Image Computing and Computer-Assisted Intervention --
  MICCAI 2015}, pages 234--241, Cham, 2015. Springer International Publishing.

\bibitem{SCHMIDHUBER2015}
J.~Schmidhuber.
\newblock Deep learning in neural networks: An overview.
\newblock {\em Neural Networks}, 61:85--117, 2015.

\bibitem{SutskeverNIPS2014}
I.~Sutskever, O.~Vinyals, and Q.~V. Le.
\newblock Sequence to sequence learning with neural networks.
\newblock In Z.~Ghahramani, M.~Welling, C.~Cortes, N.~D. Lawrence, and K.~Q.
  Weinberger, editors, {\em Advances in Neural Information Processing Systems
  27}, pages 3104--3112. Curran Associates, Inc., 2014.

\bibitem{Ulyanov2018}
D.~Ulyanov, A.~Vedaldi, and V.~Lempitsky.
\newblock Deep image prior.
\newblock {\em arXiv:1711.10925}, 2018.

\bibitem{WangChung2018}
Y.~Wang, C.~W. Siu, E.~T. Chung, Y.~Efendiev, and M.~Wang.
\newblock Deep multiscale model learning.
\newblock {\em arXiv preprint arXiv:1806.04830}, 2018.

\bibitem{yavneh2006multilevel}
I.~Yavneh and G.~Dardyk.
\newblock A multilevel nonlinear method.
\newblock {\em SIAM journal on scientific computing}, 28(1):24--46, 2006.

\end{thebibliography}
\end{document}